\documentclass[10pt]{article}

\usepackage{tikz,color,graphics,stmaryrd,blkarray,hyperref,amsmath,esint, amssymb}
\usepackage{Commons}
\usetikzlibrary{matrix}

\title{An introduction to random matrix theory}
\date{\vspace{-5ex}}

\newtheorem{theorem}{Theorem}[section] 
\newtheorem{lemma}[theorem]{Lemma}     
\newtheorem{corollary}[theorem]{Corollary}

\newtheorem{definition}[theorem]{Definition}

\renewcommand\newpage{}

\begin{document}

\maketitle

\begin{center}
Ga\"etan Borot\footnote{Max-Planck Institut f\"ur Mathematik, Vivatsgasse 7, 53111 Bonn, Germany. \\  \texttt{gborot@mpim-bonn.mpg.de}}
\end{center}

\setcounter{tocdepth}{2}
\tableofcontents

\section{Preface}
These are lectures notes for a 4h30 mini-course held in Ulaanbaatar, National University of Mongolia, August 5-7th 2015, at the summer school \textbf{Stochastic Processes and Applications, Mongolia}. The aim was to present an introduction to basic results of random matrix theory and some of its motivations, targeted to a large panel of students coming from statistics, finance, etc. Only a small background in probability is required (Mongolian students had a 1.5 month crash course on measure theory before the summer school). A few references to support -- or go further than -- the course:

\vspace{0.2cm}

\begin{itemize}
\item \emph{High Dimensional Statistical Inference and Random Matrices}, I.~Johnstone, Proceedings of the ICM, Madrid, Spain, (2006), math.ST/0611589. A short review of the application of random matrix theory results to statistics.
\item \emph{Theory of finance risks: from statistical physics to risk management}, J.~P.~Bou\-ch\-aud and M.~Potters, CUP (2000). A book explaining how ideas coming from statistical physics (and for a small part, of random matrices) can be applied to finance, by two pioneers. J.~P.~Bouchaud founded a hedge fund (Capital Fund Management), which conduct investment using those ideas, as well as pure research.
\item \emph{Population structure and eigenanalysis}, N.~Patterson, A.~L.~Preis and D.~Reich, PLoS Genetics \textbf{2} 12 (2006). Research  discussing the methodology of PCA, and proposing statistical tests based on Tracy-Widom distributions, with applications to population genetics in view.
\item \emph{Random matrices}, M.~L.~Mehta, 3rd edition, Elsevier (2004). Written by a pioneer of random matrix theory. Accessible at master level, rather focused on calculations and results for exactly solvable models, including Gaussian ensembles. A good reference to browse for results.
\end{itemize}

\noindent \textbf{Acknowledgments} I thank Carina Geldhauser, Andreas Kyprianou, Tsogzolmaa Saizmaa and the local organizers in Mongolia to have arranged this event, as well as the DAAD, the University of Augsburg and Lisa Beck for funding.

\newpage

\section{Motivations from statistics for data in high dimensions}

Collecting a huge amount of data has been facilitated by the development of computer sciences. It is then critical to have tools to analyze these data. Imagine that for each sample one has collected information represented by a point in $\mathbb{R}^N$. With $N$ large, this is certainly too much information for our brain to process. One would like to know if some relevant patterns can be identified, that would explain most of the scattering of the data by restricting to a well-chosen $k$-dimensional plane in $\mathbb{R}^N$, for $k = 1,2,3$ etc. This problem is posed for instance in archeology, in biology and genetics, in economics and finance, in linguistics, etc. Let us give some examples.

\subsection{Latent semantics}

Imagine we have documents $i \in \{1,\ldots,n + 1\}$, that we would like to group by similarity of topic. 
One strategy is to spot certain words $j \in \{1,\ldots,p\}$ in these documents, and compute the frequency $f_{ij}$ -- this can be automatized efficiently -- of occurrence of the word $j$ in document $i$. We then form the $n \times p$ matrix $X$ whose $(i,j)$-th entry is:
\beq
\label{xentry} x_{ij} = f_{ij} - \frac{1}{n + 1} \sum_{k = 1}^{n + 1} f_{kj}\,.
\eeq
Since we subtracted the mean frequency, the data $x_{n + 1,j} = -\frac{1}{n + 1}\sum_{k = 1}^{n} f_{kj}$ is determined by the $x_{ij}$ for $i \leq n$, so it is enough to consider a $n \times p$ matrix.

\begin{figure}[h!]
\begin{center}
\includegraphics[width=0.6\textwidth]{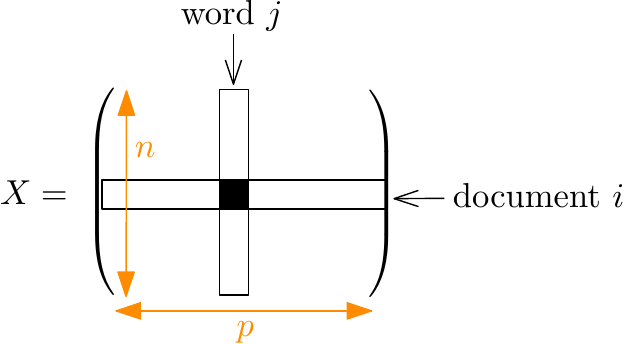}
\end{center}
\end{figure}

Let us consider the covariance matrix $M = p^{-1} XX^{T}$ ($X^{T}$ is the transpose of the matrix $X$). $M$ is a symmetric matrix of size $n \times n$, with entries:
$$
M_{ik} = \frac{1}{p} \sum_{j = 1}^{p} x_{ij}x_{kj}\,.
$$
$M_{ik}$ is large when, there are many words $j \in \{1,\ldots,p\}$ whose frequency is above the mean both in document $i$ and $k$, or below the mean both in $i$ and $k$. So, $M_{ik}$ can be considered as a measure of the correlation between the documents. For instance, if two documents both contain many "horse" and "ger", but very few "kangaroo" and "bush", the corresponding entry in the matrix $M$ will at least be made of 4 large positive terms. On the other hand, there might be many words -- for instance "river", "road", "car", "bird" -- whose frequency is close to what can be expected in an arbitrarily chosen document (clearly, one should not choose such generic words, unless one expects them for some reason to be able to differentiate the documents one wants to analyze) ; and some other words -- "tea", "cheese", "mountain" -- may sometimes appear in excess, or not very frequently, so that the sign of $x_{ij}x_{kj}$ is sometimes positive and negative without a clear trend : in these two cases, the total contribution of these words to $M_{ik}$ will be small in absolute value.

Instead of trying to group documents one by one when we notice a strong correlation -- as one can read from the large matrix $M$ -- one introduces the notion of \textbf{weighted document}, i.e. the assignment of real numbers $w_i$ to each document $i$. They can be collected in a column vector $W = (w_i)_{1 \leq i \leq n}$. Actually, only the relative weight of $i$ and $j$ matters: for any $\lambda > 0$, $W$ and $\lambda W$ represent the same weighted document. A way to fix this ambiguity is to restrict ourselves to vectors $W$ with unit euclidean norm:
$$
W^{T}W = \sum_{i = 1}^n w_i^2 = 1\,.
$$
Then, only $W$ and $-W$ represent the same weighted document. Let us try to find the weighted document $W$ that would display the strongest correlation, i.e. we want to maximize:
$$
W^{T}MW = \sum_{i,k = 1}^n w_iw_k\,M_{ik}\,.
$$
among vectors of unit norm. The answer is that $W$ should be an eigenvector\footnote{Remember that a symmetric matrix of size $n \times n$ with real-valued entries has exactly $n$ real eigenvalues, counted with multiplicity. In particular, there is a maximum eigenvalue. Besides, we are here looking at a covariance matrix, so its eigenvalues are non-negative.} of $M$ with maximum eigenvalue:
$$
MW^{(1)} = \lambda_{1}W^{(1)}\,.
$$
If $W^{(1)}_{i}$ and $W^{(1)}_{j}$ are both large and positive -- or both large and negative --  we can interpret documents $i$ and $j$ as being "similar" according to the strongest pattern that has been found in the data. If $W^{(1)}_{i}$ is close to $0$, it means that the document $i$ does not really participate to this strongest pattern.

We could also have a look at the second, the third, etc. strongest patterns, i.e. consider the eigenvectors $W^{(a)}$ for the $a$-th eigenvalue, sorted in decreasing order $\lambda_{1} \geq \lambda_{2} \geq \cdots \geq \lambda_{n}$. Unless the matrix $M$ enjoys for a special reason extra symmetries on top of $M^{T} = M$, the $n$ eigenvalues computed from the numerical data of $M$ will most likely be distinct, so there is for each $\lambda^{(a)}$ a unique (up to overall sign) eigenvector $W^{(a)}$. Let $E_a = {\rm span}(W^{(a)})$ be the eigenspace for $\lambda_{a}$. This method provides a decomposition of the space of weighted documents $\mathbb{R}^n$  into subspaces $E_{1}$, $E_{1} \oplus E_{2}$, $E_{1} \oplus E_{2} \oplus E_{3}$, \ldots of dimension $1,2,3$,\ldots In other words, it achieves the task of identifying some low dimensional subspaces $\bigoplus_{\lambda > x} E_{\lambda}$ in the high dimensional $\mathbb{R}^n$, and the threshold $x$ and dimension gives an indication of the relevance of the pattern that are identified in this way. This method is called \textbf{Principal Component Analysis (PCA)}, and was introduced in statistics by Pearson in 1901 \cite{Pearson} and Hotelling in 1931 \cite{Hotelling}.

To present PCA results, it is customary to draw in the $2$-dimensional plane a point $p_i = (x_i,y_i)$ with coordinates  $x_i = W^{(1)}_{i}$ and $y_i = W^{(2)}_{i}$ for each $i \in \{1,\ldots,n\}$. The documents that appear in the same region are then interpreted as "similar" (see Figure~\ref{Patterson}).

\subsection{Population genetics}

\begin{figure}[h!]
\begin{center}
\includegraphics[width=\textwidth]{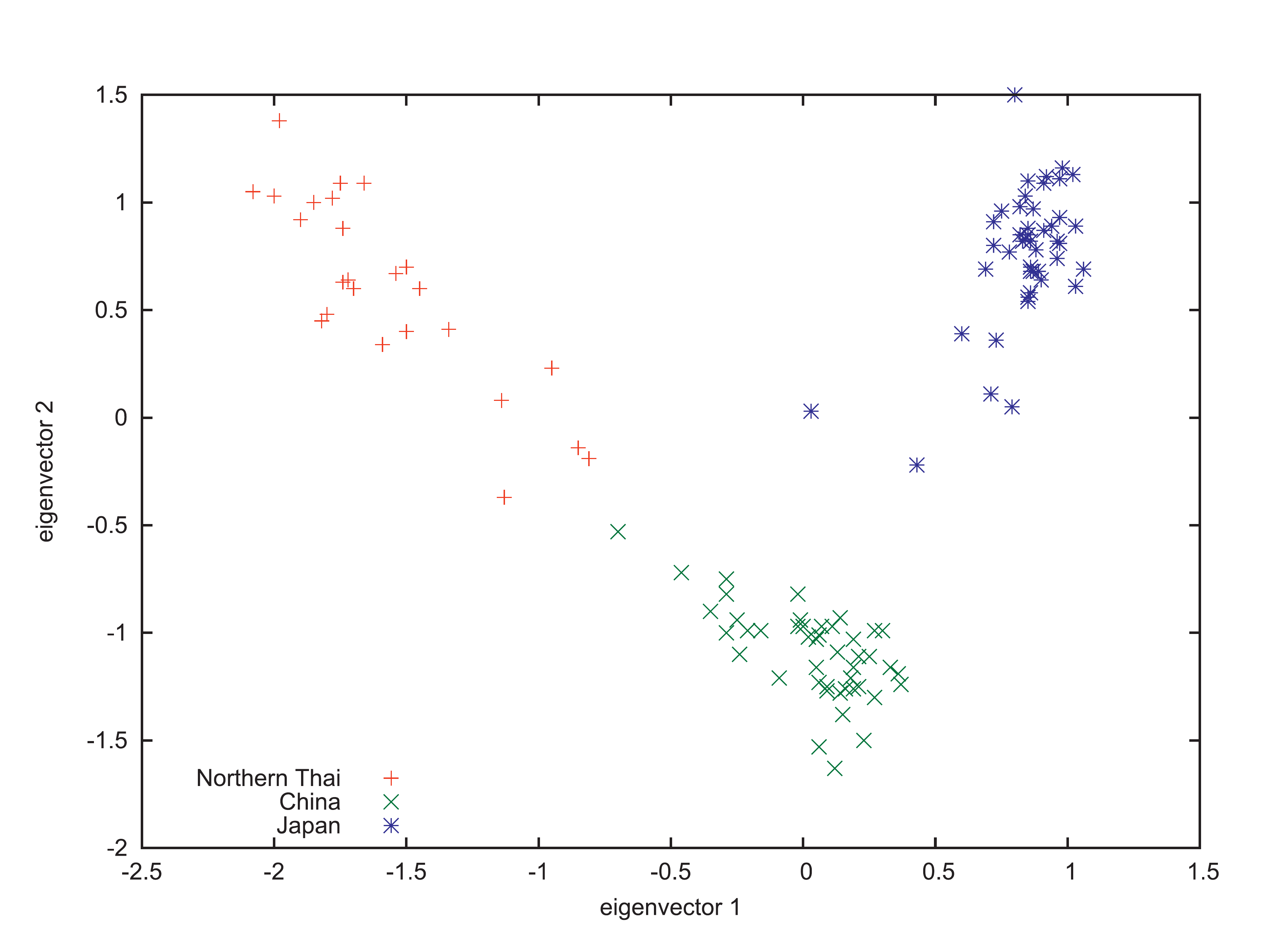}
\caption{\label{Patterson} PCA analysis of genetic data of individuals from $3$ East Asian populations, based on the International Haplotype Map, and concerning $p = 40560$ SNPs. SNP stands for Single Nucleotide Polymorphism: genes come in several versions, which often differ by the nature of the nucleotide (A, C, G or T) present in a few specific positions in the gene. Up to a correction factor, $f_{i,j}$ in \eqref{xentry} measures the frequency of a given allele (=version of a gene) $j$ carried by an individual $i$, and therefore takes values $0$, $1$ or $2$ (this last case means that the two chromosomes carry the same allele). Reprinted from \emph{Population structure and eigenanalysis}, N. Patterson, A.L. Preis and D. Reich, PLoS Genetics \textbf{2} 12 (2006).}
\end{center}
\end{figure}

If one replaces "document" by "individual", and "word" by allele (i.e. version) of a gene, the same strategy allows to study the genetic proximity of various populations, and maybe gain some insight into the history of population mixtures. Figure~\ref{Patterson} is drawn from such an example. 

\subsection{A remark}

From the matrix $X$, one could also build a $p \times p$ covariance matrix, whose lines and columns are indexed by words (or genes):
$$
\tilde{M} = n^{-1}\,X^{T}X\,.
$$
Its PCA analysis is useful for factor analysis, i.e. to study what are the most prominent reasons of similarity among the documents (or individuals).

\subsection{A word of caution}

As in any statistical analysis, care should be taken before drawing any conclusion of a cloud of points. PCA has a wide scope of applications in various disciplines, and as a result of its popularity, some research works which use PCA are not free of basic methodology errors. For instance, the most obvious fact is that points gathered near $(0,0)$ do not represent any information, except that the patterns identified do not allow to distinguish those documents. Another common mistake is to display, say  $W^{(3)}$ in abscissa and (to exaggerate) $W^{(18)}$, without questioning the relevance of the eigenvector for the $18$-th eigenvalue. It is totally possible that a very small number -- like $0$, $1$, $2$, ... -- of eigenvectors are actually relevant, the other being not distinguishable from those of a matrix with random entries.

\subsection{The use of random matrix theory}

Random matrix theory provides statistical tests for the relevance of PCA results, as follows. One chooses a \textbf{null model}, which in the previous examples would be an ensemble of symmetric random matrices $M^{{\rm null}}$. The idea behind the choice of the null model is that sampling $M^{{\rm null}}$ in this random ensemble will produce data that "contain no information" compared to the type of information we would like to identify in genuine data. Imagine that one has computed the probability $p_{A}^{{\rm null}}$ of various events $A$ concerning the eigenvalues or the eigenvectors of a matrix $M^{{\rm null}}$ drawn from the null model. If one observes the event $A$ in the genuine data one is analyzing, we say that the null model can be rejected with confidence $1 - p_{A}^{{\rm null}}$.

To this end, for various random ensembles of matrices (that one could take as null models):
\begin{itemize}
\item[$\bullet$] we need to know the distribution of eigenvalues, especially in the limit of matrices of large size ;
\item[$\bullet$] we are especially interested in extreme (maximal or minimal) eigenvalues ;
\item[$\bullet$] and we would like to understand whether these distributions are very sensitive or not to the choice of the null model, i.e. what happens to the spectrum if we do small perturbations of our random matrix.
\end{itemize}
These questions are a priori non obvious to answer, and represent typical interests in random matrix theory. 

\newpage

\section{General principles}

We shall introduce in Section~\ref{WishartS} and \ref{GaussianS} two ensembles of random matrices, but before that, let us pose the problem in mathematical terms.

\subsection{Definition and tools}

We say that a $n \times n$ matrix $M$ is symmetric if $M_{ij}$ is real and $M_{ij} = M_{ji}$, and that is hermitian is $M_{ij}$ is complex and $M_{ij} = M_{ji}^*$ where the $^*$ stands for complex conjugate. We denote:
\beq
\mathcal{S}_{n} = \big\{n \times n\,\,{\rm symmetric}\,\,{\rm matrices}\big\},\qquad \mathcal{H}_{n} = \big\{n \times n\,\,{\rm hermitian}\,\,{\rm matrices}\big\}
\eeq
and we note that $\mathcal{S}_{n} \subseteq \mathcal{H}_{n}$. The Lebesgue measure on $\mathcal{S}_n$ is by definition the product of the Lebesgue measures on the linearly independent entries of $M$:
$$
\dd M = \prod_{1 \leq i < j \leq n}\!\!\!\! \dd M_{ij}\,\prod_{i = 1}^n \dd M_{ii}\,.
$$
Similarly on $\mathcal{H}_{n}$:
$$
\dd M = \prod_{1 \leq i < j \leq n} \!\!\!\! \dd(\mathrm{Re}\,M_{ij})\,\dd(\mathrm{Im}\,M_{ij})\,\prod_{i = 1}^n \dd M_{ii}\,.
$$

A matrix $M \in \mathcal{H}_{n}$ has exactly $n$ real eigenvalues, that we write in decreasing order:
$$
\lambda_{1}^{(M)} \geq \lambda_{2}^{(M)} \geq \cdots \geq \lambda_{n}^{(M)}\,.
$$
The spectral measure is the probability measure:
$$
L^{(M)} = \frac{1}{n} \sum_{i = 1}^n \delta_{\lambda_i^{(M)}}\,.
$$
consisting of a Dirac mass $1/n$ on each eigenvalue. This is a convenient way to collect information on the spectrum of $M$, since for any continuous function $f$, we can write:
$$
\sum_{i = 1}^n f(\lambda_i^{(M)}) = \int f(x)\,\dd L^{(M)}(x)\,.
$$

We state without proof the \textbf{Hoffman-Wielandt inequality}:
$$
\forall A,B \in \mathcal{H}_{n},\qquad \sum_{i = 1}^n (\lambda_i^{A} - \lambda_i^{B})^2 \leq \mathrm{Tr}\,(A - B)^2\,.
$$
The right-hand side can be written in several forms:
$$
\mathrm{Tr}\,M^2 = \sum_{i = 1}^n (\lambda_i^{(M)})^2 = \sum_{i,j = 1}^n M_{ij}M_{ji} = \sum_{i,j = 1}^n |M_{ij}|^2\,.
$$
We remark that, since $A$ and $B$ a priori do not commute, $\lambda_i^{A} - \lambda_i^{B}$ is not in general an eigenvalue of $A - B$. This inequality is pretty useful. For instance, it tells us that the vector of eigenvalues $(\lambda_1^{(M)},\ldots,\lambda_n^{(M)})$ is a Lipschitz -- and a fortiori, continuous\footnote{Another way to prove this is to remark that the eigenvalues of $M$ are the roots of the characteristic polynomial $\det(z - M)$. The coefficients of this polynomial of $z$ are polynomial functions of the entries of $M$, thus continuous, and it is a standard result of complex analysis that the roots of a polynomial are continuous functions of the coefficients.} -- function of the entries of $M$.

\subsection{Random matrices, topology, convergence}

By convention, any topological space is equipped with the $\sigma$-algebra generated by its open sets -- the so-called Borel $\sigma$-algebra.

A random matrix of size $n$ is a random variable $M_n$ with values in $\mathcal{H}_{n}$, i.e. a measurable function from a set $\Omega$ to $\mathcal{H}_{n}$. Since eigenvalues are continuous functions of the entries, the $\lambda_i^{(M_n)}$ are also random variables, i.e. measurable functions from $\Omega$ to $\mathbb{R}$. The random probability measure $L^{(M_n)}$ is called the \textbf{empirical (spectral) measure}. At this point we need to specify the topology we choose on the set $\mathcal{M}^1(\mathbb{R})$ of probability measures on $\mathbb{R}$. We shall be concerned with two choices: the \textbf{weak topology} and the \textbf{vague topology}. For the weak topology, $\mathcal{M}^1(\mathbb{R})$ is a Polish space ; as a consequence (or as a fact for those who are not familiar with topology), it is enough to declare what does it mean for a sequence $(\mu_n)_{n}$ of probability measures to converge to a probability measure $\mu$ in this topology:
$$
\mu_n \mathop{\longrightarrow}^{{\rm weak}}_{n \infty} \mu \quad\qquad \Longleftrightarrow \quad\qquad \forall f \in \mathcal{C}^0_{b},\qquad \lim_{n \rightarrow \infty} \int f\dd\mu_n = \int f\dd\mu\,,
$$
where $\mathcal{C}^0_{b}$ is the set of continuous bounded functions from $\mathbb{R}$ to $\mathbb{R}$. For the vague topology, the convergence of sequences is nearly the same:
$$
\mu_n \mathop{\longrightarrow}^{{\rm vague}}_{n\infty} \mu\quad \qquad \Longleftrightarrow\quad\qquad \forall f \in \mathcal{C}^0_{c},\qquad \lim_{n \rightarrow \infty} \int f\dd\mu_n = \int f\dd\mu\,,
$$
where $\mathcal{C}^0_{c}$ is the set of continuous functions with compact support. Therefore, convergence for the weak topology implies convergence for the vague topology, but the converse may not hold. Now, if we equip $\mathcal{M}^1(\mathbb{R})$ is equipped with the Borel $\sigma$-algebra of any of these topologies, the empirical measure $L^{(M_n)}$ is a (probability measure)-valued random variable, i.e. a measurable function $\Omega \rightarrow \mathcal{M}^1(\mathbb{R})$. 

Usually, we are dealing with an ensemble of random matrices for each $n$, and want to study the spectrum when $n \rightarrow \infty$. We should distinguish:
\begin{itemize}
\item[$\bullet$] \textbf{global information}, which involve the macroscopic behavior of eigenvalues. For instance, we ask about the convergence of $L^{(M_n)}$ -- as a random variable -- towards a deterministic limit, its fluctuations, etc.
\item[$\bullet$] and \textbf{local information}, which concern only $O(1)$ eigenvalues. For instance, we ask about the convergence of the maximal eigenvalue $\lambda_1^{(M_n)}$, its fluctuations, etc.
\end{itemize}
We remind that, if $(X_n)_{n}$ is a sequence of random variables with values in $\mathcal{X}$, there are several (non-equivalent) notions of convergence to another $\mathcal{X}$-valued random variable $X$. The three main ones we shall use are \textbf{almost sure} convergence, convergence in \textbf{probability} and for $\mathcal{X} = \mathbb{R}$, convergence \textbf{in law}. The definitions are "$(X_n)_n$ converges to $X$ \ldots"
\begin{itemize}
\item[$\bullet$] almost surely, if $\mathbb{P}\big[\lim_{n \rightarrow \infty} X_n = X\big] = 1$.
\item[$\bullet$] in probability, if for any $\epsilon > 0$, $\lim_{n \rightarrow \infty} \mathbb{P}\big[|X_n - X| > \epsilon] = 0$.
\item[$\bullet$] in law, if for any $x \in \mathbb{R}$ at which $\mathbb{P}[X \leq x]$ is continuous, $$\lim_{n \rightarrow \infty} \mathbb{P}[X_n \leq x] = \mathbb{P}[X \leq x].$$
\end{itemize}
We remind that almost sure convergence implies convergence in probability, and the latter implies convergence in law, but the converse in general do not hold.

Even if the entries $M_{ij}$ are independent random variables, the eigenvalues depend in a non-linear way of all the entries, and therefore are strongly correlated. For this reason, the limit distributions of the spectrum in the limit $n \rightarrow \infty$ are in general very different than the limit distributions one can find in the theory of independent random variables\footnote{For independent identically distributed random variables, we have the law of large numbers and the central limit theorem for the sum, and we also know that the possible limit distributions for the maximum of a sequence of i.i.d. are the Gumbel law (e.g. for variables whose distribution decays exponentially), the Fr\'echet law (e.g. for heavy tailed distributions) and the Weibull law (e.g. for bounded random variables).}. We will see a few of these new limit laws in the lectures. It turns out these laws enjoy some universality, and the results of random matrix theory have found applications way beyond statistics, e.g. in biology and the study of ARN folding, in number theory, in nuclear physics, statistical physics and string theory, etc.

\subsection{Qualitative remarks}
\label{Remrk}
\subsubsection{Size of the spectrum}

Imagine that one fills a hermitian matrix $M_n$ of size $n$ with entries of size $O(1)$. How large (as a function of $n$) in absolute value can we expect the eigenvalues to be? We have:
$$
\mathrm{Tr}\,M_n^2 = \sum_{i,j = 1}^n \big|[M_n]_{ij}\big|^2 = \sum_{i = 1}^n \big[\lambda^{(M_n)}_i\big]^2\,.
$$
This quantity is of order $n^2$, since in the first expression it is written as a sum of $n^2$ terms of order $1$. Then, from the second expression we deduce roughly that the eigenvalues should be order $\sqrt{n}$. In other words, if we fill a matrix $M_n$ of size $n$ with entries of size $O(n^{-1/2})$ -- or equivalently with random variables having variance of order of magnitude $1/n$ -- we can expect the spectrum to remain bounded when $n \rightarrow \infty$. This non-rigorous argument serves as an explanation of the scalings chosen in the forthcoming definitions.

\subsubsection{Stability under perturbations}

Let $M_n$ be a random matrix of size $n$, and assume that when $n \rightarrow \infty$, $L^{(M_n)}$ converges to a deterministic limit $\mu$ in probability for the vague topology, i.e. for any $\epsilon > 0$ and $f \in \mathcal{C}_{c}^0$,
\beq
\label{proabc}\lim_{n \rightarrow \infty} \mathbb{P}\Big[\Big|\int f(x)\,\dd(\mu_n - \mu)(x)\Big| > \epsilon\Big] = 0\,.
\eeq
Then, let $\Delta_n$ be another random matrix of size $n$.
\begin{lemma}
If $\lim_{n \rightarrow \infty} n^{-1}\mathbb{E}[\mathrm{Tr}\,\Delta_n^2] = 0$, then $L^{(M_n + \Delta_n)}$ converges to $\mu$ in probability, for the vague topology.
\end{lemma}
\noindent \textbf{Proof.}  Any continuous $f$ with compact support can be approximated for the sup norm by a polynomial (Stone-Weierstra\ss{} theorem), in particular by a Lipschitz function. Therefore, it is enough to prove that \eqref{proabc} holds for $\mu_n = L^{(M_n + \Delta_n)}$ for any $\epsilon > 0$ and $f$ Lipschitz. Let us denote $k$ its Lipschitz constant. We have:
\bea
\Big| \int f(x)\dd(L^{(M_n + \Delta_n)} - \dd L^{(M_n)})(x)\Big| & = & \frac{1}{n} \Big| \sum_{i = 1}^n f(\lambda_i^{(M_n + \Delta_n)}) - f(\lambda_i^{(M_n)})\Big| \nonumber \\
& \leq & \frac{1}{n} \sum_{i = 1}^n k \big|\lambda_i^{(M_n + \Delta_n)} - \lambda_i^{(M_n)}\big| \nonumber \\
& \leq & \frac{k}{\sqrt{n}}\Big(\sum_{i = 1}^N (\lambda_i^{(M_n + \Delta_n)} - \lambda_i^{(M_n)})^2\Big)^{1/2} \nonumber \\
& \leq & \frac{k}{\sqrt{n}}\big(\mathrm{Tr}\,\Delta_n^2\big)^{1/2}\,, \nonumber
\eea
where we have used Cauchy-Schwarz inequality, and the Hoffman-Wielandt inequality. Then, for any fixed $\epsilon > 0$, with Markov inequality:
$$
\mathbb{P}\Big[\Big|\int f(x)\,\dd(L^{(M_n + \Delta_n)} - L^{(M_n)})(x)\Big| > \epsilon\Big] \leq \frac{k^2\,\mathbb{E}[\mathrm{Tr}\,\Delta_n^2]}{n\epsilon^2}\,,
$$
and under the assumption of the lemma, the right-hand side converges to $0$. Since we already had \eqref{proabc} for $\mu_n = L^{(M_n)}$, we have proved the desired result. \hfill $\Box$

\vspace{0.2cm}

As we have seen before, it is natural to consider matrices $M_n$ whose entries have variance bounded by $C/n$. In that case, according to this lemma, we could make $o(n^2)$ entries deterministic -- by choosing $[\Delta_n]_{ij} = \mathbb{E}[[M_n]_{ij}] - [M_n]_{ij}$ for the selected entries -- without affecting the convergence of the empirical measure to the limit $\mu$. This lemma indicates that small perturbations of a random matrix do not affect global properties of the spectrum.

There is no such general rule for local properties (such as the position of the maximum eigenvalue): we will see examples showing that sometimes they are preserved under small perturbations, and sometimes they are dramatically affected.

\newpage

\section{Wishart matrices}
\label{WishartS}
\subsection{Definition}

A \textbf{real Wishart matrix} is a random symmetric matrix $M$ of the form:
$$
M = n^{-1}\,X^{T}X\,,
$$
where $X$  is random  matrix of size $n \times p$ such that:
\begin{itemize}
\item[$\bullet$] $(X_{ij})_{1 \leq i \leq n}$ are independent samples of a real-valued random variable $\mathcal{X}_j$ ; 
\item[$\bullet$] $(\mathcal{X}_{1},\ldots,\mathcal{X}_{p})$ is a Gaussian vector with given covariance $K \in \mathcal{S}_{p}$
\end{itemize}
In other words, the joint probability density function (= p.d.f.) of the entries of $X$ is:
$$
c_{np}(K)\,\exp\Big(-\frac{1}{2} \sum_{i,i' = 1}^n \sum_{j,j' = 1}^p X_{ij}X_{i'j'} K^{-1}_{jj'}\Big) = c_{np}(K)\,\exp\Big(-\frac{1}{2} \mathrm{Tr}\,X^{T} K^{-1} X\Big)\,.
$$
$c_{np}(K)$ is a normalization constant. All the normalization constants that will appear in these lectures can be explicitly computed, but we will not care about them. The matrix $M$ is of size $p \times p$, and $n$ is called the number of degrees of freedom. The parameter:
$$
\gamma = n/p
$$
will play an important role. The ensemble of real Wishart matrices with a covariance $K = {\rm diag}(\sigma^2,\ldots,\sigma^2)$ is a natural choice of null model for covariance matrices in data analysis, which depends on a parameter $\sigma$. It was introduced by Wishart in 1928 \cite{Wishart}.

One can also define the ensemble of \textbf{complex Wishart matrices}. These are random hermitian matrices of the form $M = (X^{T})^*X$, where $(X_{ij})_{1 \leq i \leq n}$ are independent samples of $\mathcal{X}_{j}$ such that $(\mathcal{X}_{1},\ldots,\mathcal{X}_{p})$ is a complex Gaussian vector with given covariance $K \in \mathcal{H}_{p}$. This is one of the simplest model of complex random matrices, and the latter are relevant e.g. in telecommunications, when one studies non-ideal propagation of waves along many canals (complex numbers are used to encode simultaneously the amplitude and the phase of a wave).

\subsection{Spectral density in the large size limit}

\begin{figure}[h!]
\begin{center}
\includegraphics[width=0.8\textwidth]{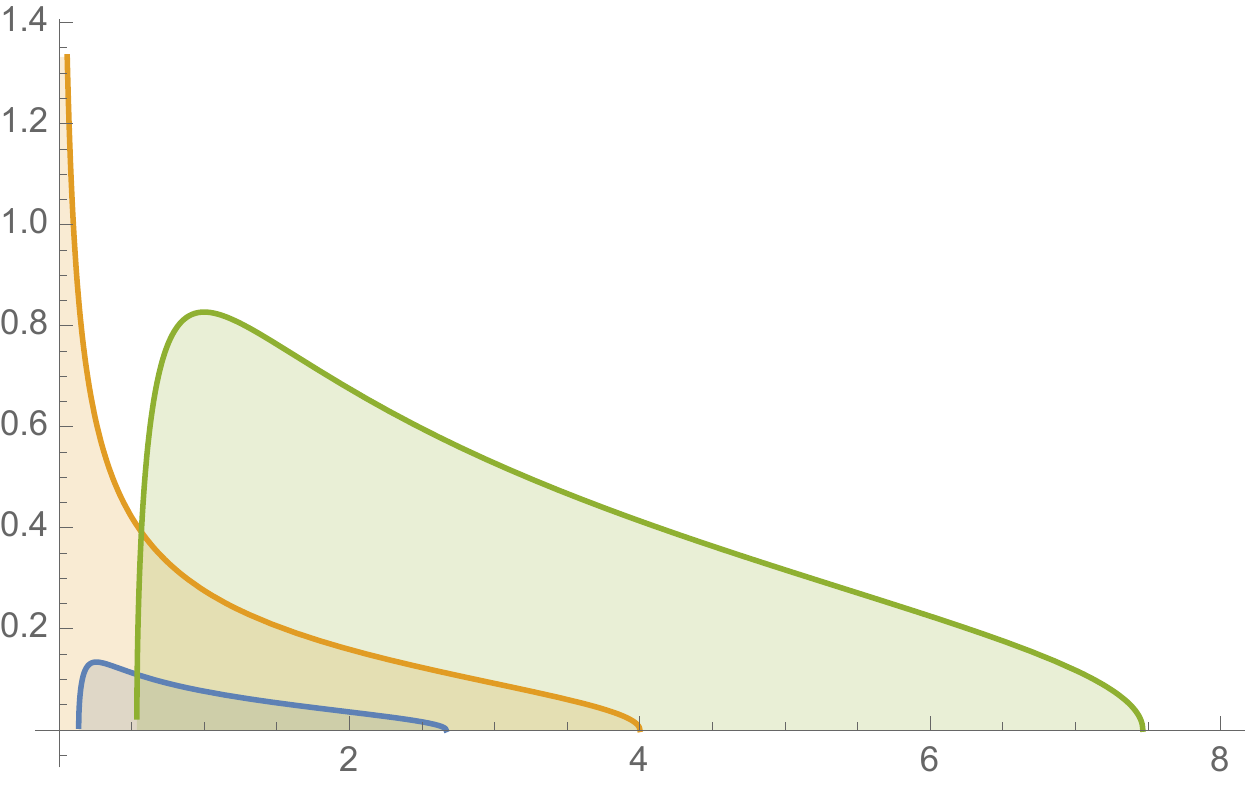}
\caption{\label{MP} Mar\v{c}enko-Pastur probability density function, for $\sigma^2 = 1$: in green $\gamma = 3$, in orange $\gamma = 1$, in blue $\gamma = 0.4$. The mass of the distribution in this last case is $0.4$, to which should be added a Dirac mass with mass $0.6$ at $0$.}
\end{center}
\end{figure}

We consider real or complex Wishart ensembles with given covariance $K = {\rm diag}(\sigma^2,\ldots,\sigma^2)$. Mar\v{c}enko and Pastur showed in 1967 \cite{MPbib} that the empirical measure $L^{(M)}$ has a deterministic limit:
\begin{theorem}
\label{MPth}
In the limit where $p,n \rightarrow \infty$ while $n/p$ converges to a fixed value $\gamma \in (0,+\infty)$, $L^{(M)}$ converges almost surely and in expectation in the weak topology, towards the probability measure (see Figure~\ref{MP}):
\beq
\label{MPlaw} \mu_{{\rm MP}} = \max(1 - \gamma,0)\delta_{0} + \frac{\gamma \sqrt{(a_+(\gamma) - x)(x - a_-(\gamma))}}{2\pi \sigma^2\,x}\,\mathbf{1}_{[a_-(\gamma),a_+(\gamma)]}\,\dd x
\eeq
where $a_{\pm}(\gamma) = \sigma^2(1 \pm \gamma^{-1/2})^2$.
\end{theorem}

We note that when $n < p$, the matrix $X^{T}X$ has rank $n < p$, and therefore has almost surely $p - n = p(1 - \gamma)$ zero eigenvalues, which explains the Dirac mass in \eqref{MPlaw} which appear for $\gamma < 1$. The mean and variance of the Mar\v{c}enko-Pastur distribution are:
\beq
\int x\,\dd\mu_{{\rm MP}}(x) = \sigma^2,\qquad \int x^2\,\dd\mu_{{\rm MP}}(x) - \Big(\int x \dd\mu_{{\rm MP}}(x)\Big)^2 = \sigma^4/\gamma\,.
\eeq 
Apart from the possible Dirac mass at $0$, the support of $\mu_{{\rm MP}}$ is spread on an interval of length $4\sigma^2\gamma^{-1/2}$ around the mean $\sigma^2$: the smaller $\gamma$ is, the broader the support becomes. On the other hand, when $\gamma \rightarrow \infty$, the support becomes localized around $\sigma^2$, i.e. we can read the variance of the Gaussian entries of $X$. For practical applications, this means that if the number of measurements $n$ is not very large compared to the number $p$ of properties we measure, the spectrum of $M$ will be spread.

Another property of $\mu_{\rm MP}$ is that, for\footnote{For $\gamma = 1$, it diverges as $x^{-1/2}$ when $x \rightarrow 0^+$.} $\gamma \neq 1$, the density of $\mu_{{\rm MP}}$ vanishes like a squareroot at the edges $a_{\pm}(\gamma)$. This behavior is frequent for the spectra of large random matrices.

\subsection{Maximum eigenvalue and fluctuations}
\label{Sectin2}
From Mar\v{c}enko-Pastur theorem, one can easily deduce that, for any $\epsilon > 0$,
$$
\mathbb{P}[\lambda_{1}^{(M)} \leq a_+(\gamma) - \epsilon] \rightarrow 0\,,
$$
and thus that $(\limsup_{n \rightarrow \infty} \lambda_{1}^{(M)})$ is almost surely larger than $a_+(\gamma)$. Indeed, let us choose an arbitrary non-negative, non-zero, continuous function $f$ with compact support included in $(a_+(\gamma) - \epsilon,+\infty)$. We can rescale $f$ to enforce $\int f(x)\dd\mu_{{\rm MP}}(x) = 1$. We then have:
\bea
\mathbb{P}\big[\lambda_{1}^{(M)} \leq a_+(\gamma) - \epsilon\big] & \leq & \mathbb{P}\Big[\int f(x)\,\dd L^{(M)}(x) = 0\Big] \nonumber \\
& \leq & \mathbb{P}\Big[\Big|\int f(x)\dd(L^{(M)} - \mu_{{\rm MP}})(x)\Big| \geq 1/2\Big]\,, \nonumber
\eea
and the latter converges to $0$ when $n,p \rightarrow \infty$ according to Theorem~\ref{MPth}. But Theorem~\ref{MPth} does not tell us whether the maximum eigenvalue $\lambda_{1}^{(M)}$ really converges to $a_+(\gamma)$ or not. The reason is easily understood: the event $\lambda_1^{(M)} \leq a_+(\gamma) - \epsilon$ actually means that all eigenvalues are smaller than $a_+(\gamma) - \epsilon$: this is a global information, hence contained in the statement of convergence of $L^{(M)}$. However, the realization of an event like $\lambda_1^{(M)} \geq a_+(\gamma) - \epsilon$ only involves a single eigenvalue, and thus more work is needed to estimate its probability. We will not say how this work is done, but the result is that there is no surprise:
\begin{theorem} \cite{Geman}
$\lambda_{1}^{(M)}$ converges almost surely to $a_+(\gamma)$.
\end{theorem}

\begin{figure}[h!]
\begin{center}
\includegraphics[width=\textwidth]{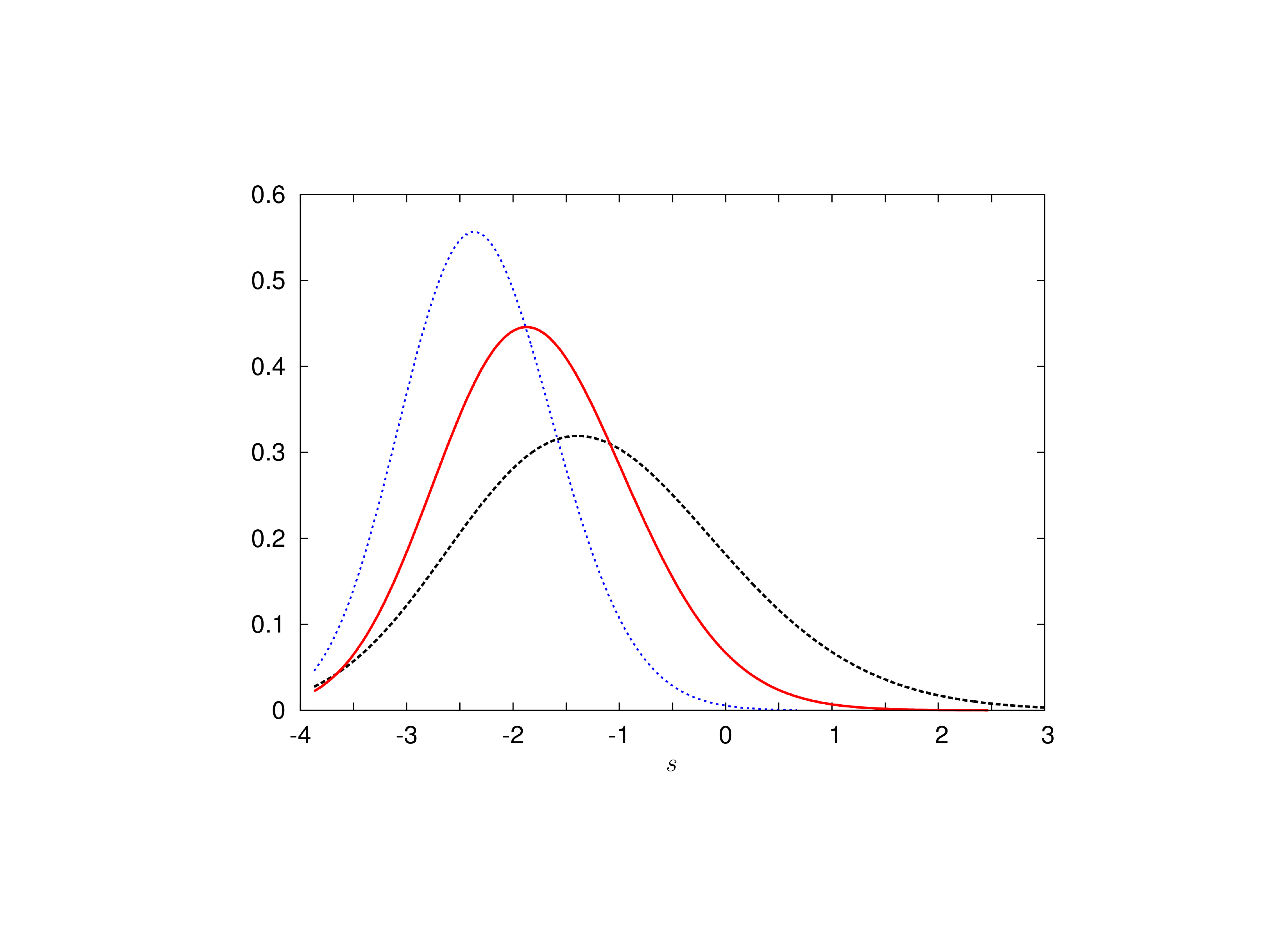}
\caption{\label{T2Wgraph} Probability density function of the Tracy-Widom law, i.e. ${\rm TW}_{\beta}(s)$, for $\beta = 1$ (GOE, in blue), $\beta = 2$ (GUE, in red), and $\beta = 4$. Graph courtesy of J.M.~St\'ephan.}
\end{center}
\end{figure}

The distribution of the fluctuations of $\lambda_1^{(M)}$ is also known. Before presenting the result, let us give a non-rigorous argument to guess the order of magnitude of these fluctuations. The guess is that, for a Wishart matrix of large size $p$, the number of eigenvalues in an interval $I_p$ whose length depend on $p$ should be well approximated by $p\mu_{{\rm MP}}[I_{p}]$. So, we guess that the fluctuations of $\lambda_1^{(M)}$ should occur in a region of width $\delta_p \rightarrow 0$ around $a_+(\gamma)$ where $\mu_{{\rm MP}}$ has mass of order $1/p$. Since $\mu_{{\rm MP}}$ vanishes like a squareroot at the edge, we have:
$$
\mu_{{\rm MP}}[a_+(\gamma) - \delta_p,a_+(\gamma)] \sim \int_{0}^{\delta_{p}} x^{1/2}\dd x = \frac{2}{3} \delta_{p}^{3/2}\,,
$$
and this gives the estimate $\delta_{p} \sim p^{-2/3}$. The following result \cite{ForresterS,JohnTW} confirms this guess:
\begin{theorem}
\label{Fluctumax}
We set $\beta = 1$ for real Wishart, and $\beta = 2$ for complex Wishart. The random variable:
$$
\gamma^{1/2}p^{2/3}\,\frac{\lambda_1^{(M)} - a_+(\gamma)}{\sigma^2(1 + \gamma^{-1/2})^{4/3}}
$$
converges in law towards a random variable $\Xi_{\beta}$ when $n,p \rightarrow \infty$ while $n/p$ converges to $\gamma \in (0,+\infty)$.
\end{theorem}
The distribution function:
$$
{\rm TW}_{\beta}(s) = \mathbb{P}[\Xi_{\beta} \leq s]
$$
is called the \textbf{Tracy-Widom law}. It is not an elementary function, but can be considered as a new special function. It is nowadays well-tabulated, hence ready for use in statistics (Figure~\ref{T2Wgraph}).
We now give one of its expression, first obtained by Tracy and Widom in 1992 for $\beta = 2$ \cite{TW92} and 1995 for $\beta = 1$ \cite{TW95}:
\bea
\label{Herm}{\rm hermitian}\qquad {\rm TW}_{2}(s) & = & \exp\Big[-\int_{s}^{\infty} \big\{q'(t) - tq^2(t) - q^4(t)\big\}\dd t\Big]\,,  \\
\label{TW1}{\rm symmetric}\qquad {\rm TW}_{1}(s) & = & \exp\Big[-\frac{1}{2}\int_{s}^{\infty} q(t)\dd t\Big]\,.
\eea
Here, $q(t)$ is the unique bounded solution to the \textbf{Painlev\'e II equation}: 
$$
q''(t) = 2q^3(t) + tq(t)
$$
satisfying the growth conditions $q(t) \sim \sqrt{-t/2}$ when $t \rightarrow -\infty$, and:
$$
q(t) \sim \frac{\exp(-\frac{2}{3}t^{3/2})}{2\sqrt{\pi}t^{1/4}},\qquad t\rightarrow +\infty\,.
$$
Existence and uniqueness of the function $q(t)$ was shown by Hastings and McLeod in 1980 \cite{HMcL}, and it bears their name. We will derive in Section~\ref{Asymmax} another expression for ${\rm TW}_{2}(s)$ in terms of a infinite size (Fredholm) determinant, which is actually the easiest way to compute numerically the Tracy-Widom law.

\subsection{Application to Markowitz portfolio optimization}
\label{MarS}
This paragraph is based on the article \emph{Random matrix theory and financial correlations}, Bouchaud, Cizeau, Laloux, Potters, Risk Magazine \textbf{12} 69 (1999), and the figures extracted from this article.

Imagine we consider investing in assets $j \in \{1,\ldots,p\}$ a fraction of money $w_j$. We would like to determine, for a fixed return $r$, the choice of portfolio $(w_1^*,\ldots,w_n^*)$ minimizing the risk. For this purpose, we only have at our disposal the observations of the price $p_{ij}$ of these assets at 
times $i \in \{1,\ldots,n\}$ in the past. We can subtract the mean price and write $p_{ij} = \overline{p}_{j} + x_{ij}$. If we had invested in the past and get our return at time $i$, we would have earned:
$$
r_i = \sum_{j = 1}^p w_j(\overline{p}_j + x_{ij})
$$
If we are ready to believe\footnote{This is highly criticizable, especially in finance. We will come back to this point.} that these observations represent well what can happen during the (future) period of our investment, we can take:
$$
r = \sum_{j = 1}^p w_j \overline{p}_{j} + \frac{1}{n} \sum_{i = 1}^n w_jx_{ij} = \overline{r} + J^{T}XW
$$
where $W$ is column vector representing the portfolio, $J$ the column vector with entries $1/n$, and $X = (x_{ij})_{ij}$ the $n \times p$ matrix collecting the observations. One can also try to evaluate the risk in investing as $W$ with the quantity:
$$
\rho = \sum_{j,j' = 1}^p w_{j}w_{j'}\Big(\frac{1}{n} \sum_{i = 1}^n x_{ij}x_{ij'}\Big) = W^{T}MW\,,
$$
where:
$$
M = n^{-1}\,X^{T}X
$$
is the empirical correlation matrix. Finding the $W^*$ that minimizes $\rho$ for a given $(r - \overline{r})$ can be done by minimizing the quantity $W^{T}MW - a J^{T}XW$ for a constant $a$ -- the Lagrange multiplier -- that we adjust so that:
$$
r - \overline{r} = J^{T}XW^*\,.
$$
Denoting $P = J^{T}X$, the result is:
\beq
\label{Marko}\rho^* = \frac{(r - \overline{r})^2}{P^{T}M^{-1}P}\,\qquad W^* = \frac{\rho^*}{r - \overline{r}}\,M^{-1}P\,.
\eeq
In particular, we see that the eigenvectors of $M$ with small eigenvalues play an important role in the evaluation of $\rho^*$ and $W^*$. This is the base of the method proposed by Markowitz in 1952 \cite{Markowitz}. One usually plots the return $r$ as a function of the estimation $\rho^*$ of the risk: the curve is called the \textbf{efficient frontier}, and in this simple model, it is a parabola.

\begin{figure}[h!]
\begin{center}
\includegraphics[width=1.2\textwidth]{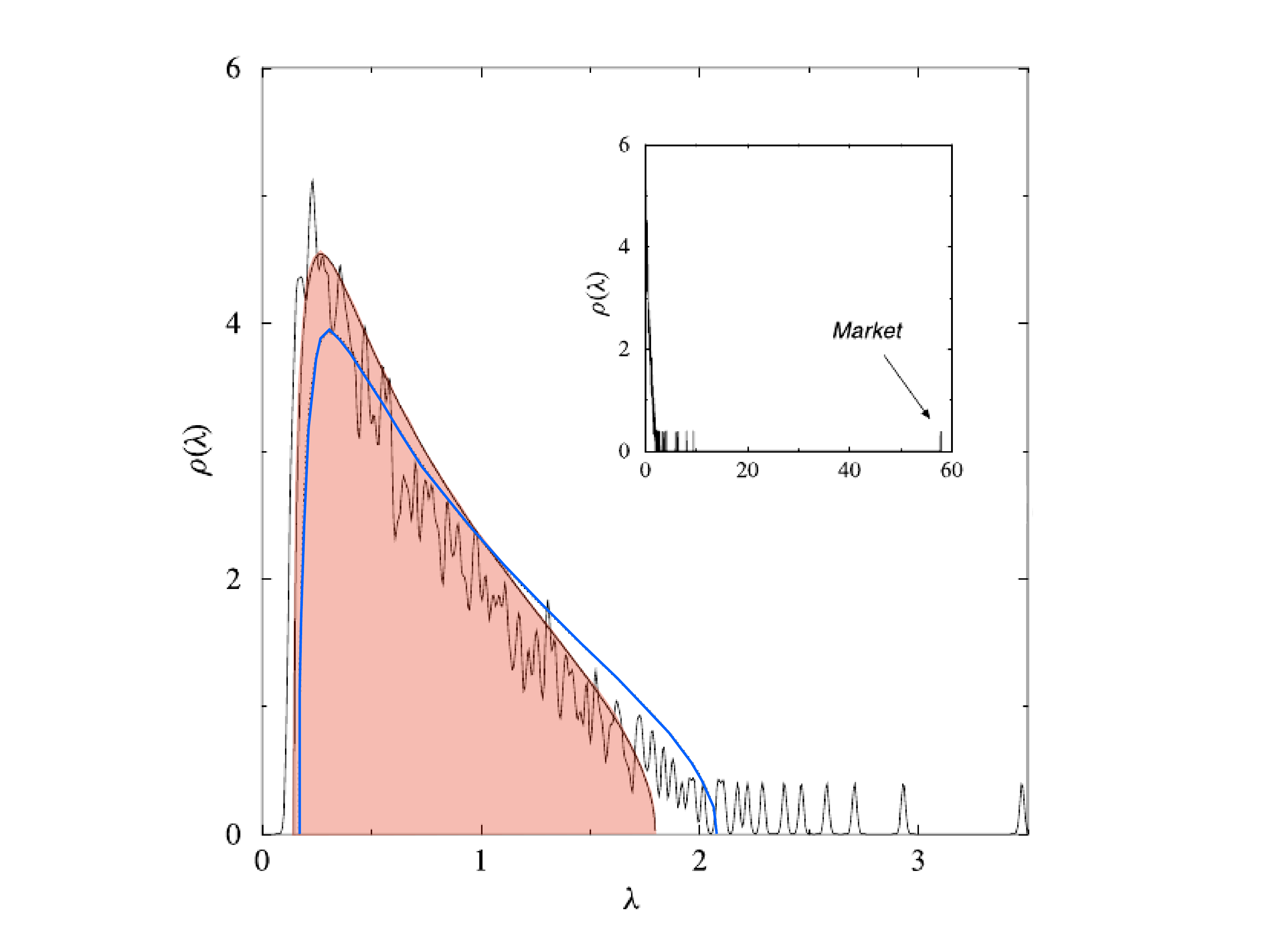}
\caption{\label{Bouchaud-1} Spectrum of an empirical $p \times p$ covariance matrix, built from the value of $p = 406$ assets from the S\&P 500, observed every day in a period of $n = 1309$ days between 1991 and 1996. One eigenvalue is much larger than the others, and correspond to the market mode, i.e. all assets increase or decrease simultaneously. The blue (resp. red) curve is the Mar\v{c}enko-Pastur (MP) spectral density for a large Wishart matrix with $\gamma = n/p$, and input covariance ${\rm diag}(\sigma^2,\ldots,\sigma^2)$ for $\sigma^2 = 0.85$ (resp. $\sigma^2 = 0.74$). This last value is the optimal fit. About $6\%$ of the eigenvalues cannot be not accounted by the MP law, and they are responsible for $1 - \sigma^2 = 26\%$ of the variance. We note that the shape of the empirical density of low eigenvalues is well reproduced by MP, so these eigenvalues (and the corresponding eigenvectors, which have the largest weight for Markowitz optimization) cannot be distinguished from noise.}
\end{center}
\end{figure}

\begin{figure}[h!]
\begin{center}
\includegraphics[width=\textwidth]{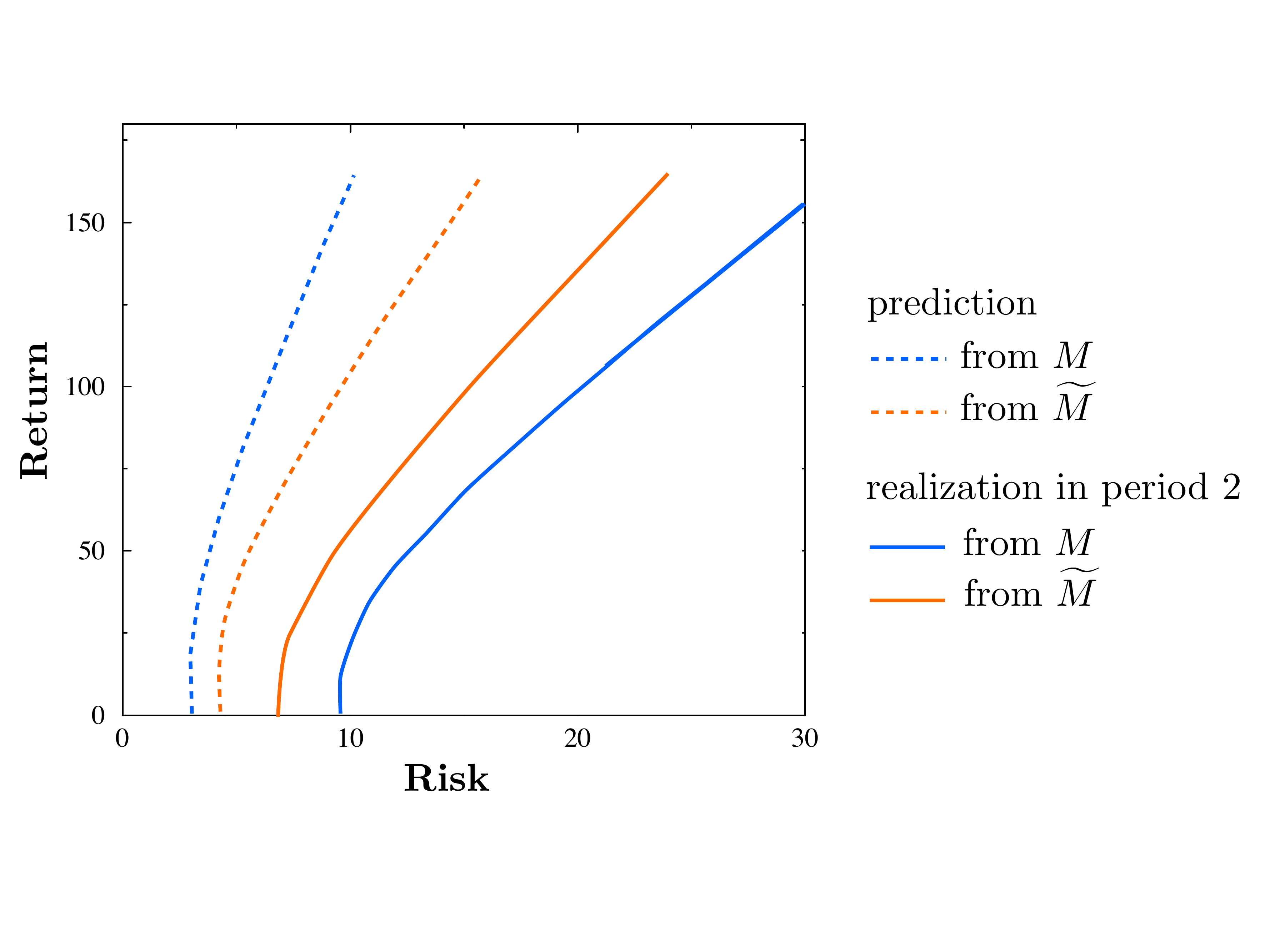}
\caption{\label{Bouchaud-2} The dashed curve is the prediction from $M$ (blue) or $\widetilde{M}$ (orange) of the effective frontier via \eqref{Marko}, constructed from the observations in a first period of time, and proposing to invest $W^*$ or $\widetilde{W}^*$. The plain curves correspond to the effective frontier measured if we really invested  $W^*$ (blue) or $\widetilde{W}^*$ (orange) in the second period of time.}
\end{center}
\end{figure}

As a matter of fact, it is hard to build an empirical covariance matrix reliable for future investments, and Markowitz theory suffers in practice from important biases. With an example drawn from genuine financial data, Bouchaud et al. pointed out that a large part -- and especially the lower part -- of the spectrum of $M$ can be fitted with a Mar\v{c}enko-Pastur distribution, hence cannot be distinguished from the null model of a large random covariance matrix (Figure~\ref{Bouchaud-1}). The effect is that the minimal risk for a given return is underestimated (Figure~\ref{Bouchaud-2}), and the guess \eqref{Marko} of the optimal portfolio does not give good results.

\begin{figure}[h!]
\begin{center}
\includegraphics[width=\textwidth]{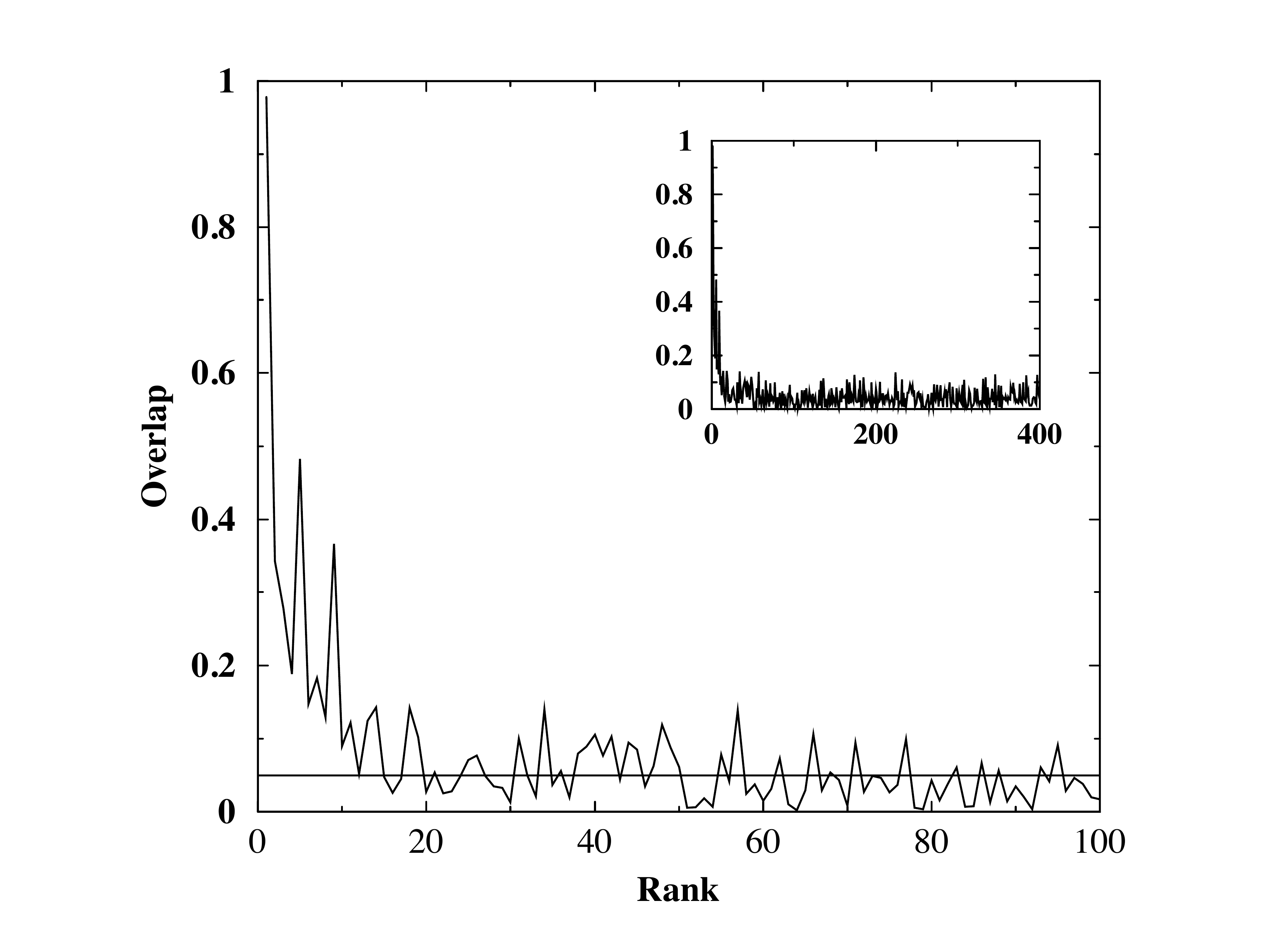}
\caption{\label{Bouchaud-3} $M^{(1)}$ and $M^{(2)}$ are empirical correlation matrices coming from observation in period $1$ and $2$. If we denote $W^{(a)}_{i}$ a unit norm eigenvector of $M^{(a)}$ for the $i$-th eigenvalue (in decreasing order) with norm $1$, the plot shows the scalar product $|W^{(1)}_{i}\cdot W^{(2)}_{i}|$ as a function of $i = 1,2,3,\ldots$ in abscissa. The horizontal line $1/\sqrt{p}$ is the typical value for the overlap of two independent random vectors with normal entries Gaussian entries.}
\end{center}
\end{figure}

The part of the spectrum undistinguishable from noise is called the \textbf{noise band}. If one makes observations of the prices and builds empirical correlation matrices over two distinct periods, one can also check that the eigenvectors for eigenvalues outside the noise band have common features -- quantitatively measured by the absolute value of their scalar product -- while the eigenvectors for eigenvalues in the noise band have nothing more in common than two random vectors (Figure~\ref{Bouchaud-3}). It supports the idea that only eigenvectors for eigenvalues outside the noise band  contain a genuine information about the long-time evolution of the market.

Although there is no ideal cure, Bouchaud et al. proposed to replace the empirical correlation matrix $M$ by $\widetilde{M}$ built as follows.
\begin{itemize}
\item[$\bullet$] Decompose $\mathbb{R}^n = E_{{\rm noise}} \oplus E$, where $E_{{\rm noise}}$ (resp. $E$) is the sum of ei\-gen\-spaces for eigenvalues in the noise band (resp. outside the noise band).
\item[$\bullet$] Replace the restriction of $M$ to $E_{{\rm noise}}$ by a multiple of the identity operator, so that the trace is preserved.
\item[$\bullet$] Use the new matrix $\widetilde{M}$ in the Markowitz optimization formulas \eqref{Marko}.
\end{itemize}
The risk is still underestimated, but to a smaller extent.

\newpage

\section{Gaussian ensembles}
\label{GaussianS}
The Gaussian ensembles are the simplest ensembles of random matrices from the computational point of view. As Wishart matrices, they come in two flavors, depending whether one considers symmetric or hermitian matrices. For a reason revealed in Section~\ref{FUm}, the symmetric case is labeled $\beta = 1$, and the hermitian case $\beta = 2$.

In the \textbf{Gaussian Orthogonal Ensemble} (GOE), we consider a symmetric random matrix $M$ of size $n \times n$, with
\beq
\label{ho}M_{ij} = \left\{\begin{array}{lll} X_{ij} & & 1 \leq i < j \leq n \\ X_{ji} & & 1 \leq j < i \leq n \\ Y_i & & 1 \leq i = j \leq n\end{array}\right.\,\,,
\eeq
where $X_{ij}$ and $Y_i$ are independent centered Gaussian random variables with:
\beq
\label{variancedef}\mathbb{E}[X_{ij}^2] = \sigma^2/n,\qquad \mathbb{E}[Y_{i}^2] = 2\sigma^2/n\,.
\eeq
We choose to scale the variance by $1/n$, so that the spectrum will remain bounded -- see Section~\ref{Remrk}. The difference of normalization between the off-diagonal and diagonal elements is motivated by observing that the resulting probability measure on the entries of $M$ is proportional to:
\beq
\label{GOElaw}\dd M\,\exp\Big[-\frac{n}{2\sigma^2}\Big(2\sum_{i = 1}^n M_{ii}^2 + \sum_{1 \leq i < j \leq n} M_{ij}^2\Big)\Big] = \dd M\,\exp\Big[-\frac{n}{2\sigma^2}\,\mathrm{Tr}\,M^2\Big]\,.
\eeq
The Lebesgue measure $\dd M$ is invariant under conjugation $M \mapsto \Omega^{-1} M \Omega$ by an orthogonal matrix $\Omega$, and so is $\mathrm{Tr}\,M^2$. Therefore, for any orthogonal matrix $\Omega$, $M$ drawn from \eqref{GOElaw} and $\Omega^{-1} M \Omega$ have the same distribution, and this explains the name GOE. This property would not be true if we had chosen the same variance in \eqref{ho} for the diagonal and off-diagonal entries.

In the \textbf{Gaussian Unitary Ensemble} (GUE), we consider a hermitian random matrix $M$ of size $n \times n$, with
$$
M_{ij} = \left\{\begin{array}{lll} X_{ij} + \sqrt{-1}\,\tilde{X}_{ij} & & 1 \leq i < j \leq n \\ X_{ji} + \sqrt{-1}\,\tilde{X}_{ji} & & 1 \leq j < i \leq n \\ Y_{ii} & & 1 \leq i = j \leq n \end{array}\right.
$$
where $X_{ij}$, $\tilde{X}_{ij}$ and $Y_{i}$ are independent centered Gaussian random variables with:
$$
\mathbb{E}[X_{ij}^2] = \mathbb{E}[\tilde{X}_{ij}^2] = \sigma^2/2n,\qquad \mathbb{E}[Y_{i}^2] = \sigma^2/n\,.
$$
The resulting probability measure on the entries of $M$ reads:
$$
\dd M\,\exp\Big[-\frac{n}{\sigma^2}\,\mathrm{Tr}\,M^2\Big]\,,
$$
and it is invariant under conjugation $M \mapsto \Omega^{-1}M\Omega$ by a unitary matrix $\Omega$.

The probability measures for the GOE and the GUE can written in a unified way:
$$
\dd M\,\exp\Big[-\frac{n\beta}{2\sigma^2} \mathrm{Tr}\,M^2\Big]\,.
$$

The results that we have seen in the case of Wishart matrices for the spectral density in the large size limit, and the location of the maximum eigenvalue and its fluctuations, have an analog for the Gaussian ensembles. Their proof in the case $\beta = 2$ (GUE) will be sketched in Section~\ref{Asym}.

\subsection{Spectral density}

Let $M_n$ be a random matrix in the GOE or the GUE. Wigner showed in 1955 \cite{Wignero} that the empirical measure $L^{(M_n)}$ converges to a deterministic limit -- although the almost sure mode of convergence was only obtained later, by large deviation techniques -- see e.g. the book \cite{AGZbook}.

\begin{theorem}
When $n \rightarrow \infty$, $L^{(M_n)}$ converges almost surely and in expectation to the probability measure (see Figure~\ref{Hermite} for a plot):
\beq
\label{SClaw}\mu_{{\rm sc}} = \frac{\sqrt{4\sigma^2 - x^2}}{2\pi\sigma^2}\,\mathbf{1}_{[-2\sigma,2\sigma]}(x)\,\dd x\,.
\eeq
\end{theorem}
$\mu_{{\rm sc}}$ is called the semi-circle law, because of the shape of its density when $\sigma = 1$. It is symmetric around $0$, and the variance is:
$$
\int_{-2\sigma}^{2\sigma} x^2\,\dd\mu_{{\rm sc}}(x) = \sigma^2\,.
$$
As in the Wishart case, we observe that the density of $\mu_{{\rm sc}}$ vanishes like a squareroot at the edges of its support. 

\subsection{Maximum eigenvalue and fluctuations}
\label{Mxun}
\begin{theorem}
\cite{MehtaBook}
When $n \rightarrow \infty$, $\lambda_{1}^{(M_n)}$ converges almost surely to $2\sigma$. Besides, we have the convergence in law:
$$
n^{2/3}\sigma^{-1}\big\{\lambda_1^{(M_n)} - 2\sigma\big\} \mathop{\longrightarrow}_{n\infty} \Xi_{\beta}\,,
$$
where $\Xi_{\beta}$ is drawn from the Tracy-Widom law with $\beta = 1$ for GOE, and $\beta = 2$ for GUE.
\end{theorem}

Comparing to the Wishart case, we remark that the global properties of the spectrum do not depend on the type -- $\beta = 1$ for symmetric, or $\beta = 2$ for hermitian -- of matrices once the ensemble is properly normalized, while the local properties (e.g. the Tracy-Widom laws) depend non-trivially on $\beta$, as one can see in Figure~\ref{T2Wgraph}.

\newpage

\section{Stieltjes transform and freeness}

\subsection{Stieltjes transform and its properties}

If $\mu$ is a probability measure on $\mathbb{R}$, its \textbf{Stieltjes transform} is the function:
\beq
\label{Stieldef} W_{\mu}(z) = \int_{\mathbb{R}} \frac{\dd\mu(x)}{z - x}\,.
\eeq
It is a holomorphic function\footnote{The support ${\rm supp}\,\mu$ is the set of all points $x \in \mathbb{R}$ such that, for any open neighborhood $U_{x}$ of $x$, $\mu[U_{x}] > 0$.} of $z \in \mathbb{C} \setminus {\rm supp}\,\,\mu$. It is an important tool because of the Stieltjes continuity theorem -- see for instance \cite{Taobook}. In its most basic form:
\begin{theorem}
Let $(\mu_n)_{n}$ be a sequence of probability measures on $\mathbb{R}$, and $\mu$ another probability measure. $\mu_n$ converges to $\mu$ for the vague topology if and only if for all $z \in \mathbb{C}\setminus\mathbb{R}$, $W_{\mu_n}(z)$ converges to $W_{\mu}(z)$.
\end{theorem}
The same theorem holds if $(\mu_n)_n$ is a sequence of random measures, by adding on both sides of the equivalence the mode of convergence "almost sure", "in probability", etc. Thus, the problem of checking the convergence of probability measures can thus be replaced with the -- usually easier -- problem of checking pointwise convergence of holomorphic functions. Let us give a few useful properties to handle the Stieltjes transform.

\noindent $\bullet$ Firstly, if $\mu$ is a measure which has moments up to order $K$, we have the asymptotic expansion:
$$
W_{\mu}(z) = \frac{1}{z} + \sum_{k = 1}^{K} \frac{m_{k}}{z^{k + 1}} + o(z^{-(K + 1)}),\qquad m_k = \int_{\mathbb{R}} x^{k}\,\dd\mu(x)
$$
valid when $|z| \rightarrow \infty$ and $z$ remains bounded away from the support (if the support is $\mathbb{R}$, that means $|{\rm Im}\,z| \geq \delta$ for some fixed $\delta > 0$). So, the moments can be read off the expansion of $W_{\mu}(z)$ at infinity.

\noindent $\bullet$ Secondly, the Stieltjes transform can be given a probabilistic interpretation. We observe that, for $y \in \mathbb{R}$ and $\eta > 0$,
$$
-\frac{1}{\pi}\,{\rm Im}\,W_{\mu}(y + {\rm i}\eta) = \int_{\mathbb{R}} \frac{\eta}{\pi}\,\frac{\dd \mu(x)}{(y - x)^2 + \eta^2}
$$ 
is the density -- expressed in the variable $y$ -- of the convolution $\mu \star C_{\eta}$ of the initial measure $\mu$ with the Cauchy measure of width $\eta$:
$$
C_{\eta} = \frac{\eta\,\dd x}{\pi(x^2 + \eta^2)}\,.
$$

\noindent $\bullet$ Thirdly, the measure $\mu$ can be retrieved from its Stieltjes transform. Indeed, if $f$ is a continuous function bounded by a constant $M > 0$, we know that:
$$
\lim_{\eta \rightarrow 0} \int_{\mathbb{R}} \frac{\eta}{\pi}\,\frac{f(y)\,\dd y}{(x - y)^2 + \eta^2} = f(x)\,,
$$
and actually the quantity inside the limit is bounded by $M$. So, by dominated convergence, we have:
\bea
\label{swesk} 
 \lim_{\eta \rightarrow 0^+} \int_{\mathbb{R}} f(y)\,\dd(\mu \star C_{b})(y) & = & \lim_{\eta \rightarrow 0^+} \int_{\mathbb{R}} \dd \mu(x)\Big(\int_{\mathbb{R}} \frac{\eta}{\pi}\,\frac{f(y)\,\dd y}{(x - y)^2 + \eta^2}\Big) \\
& = & \int_{\mathbb{R}} f(x)\,\dd\mu(x)\,. \nonumber
\eea
This means that, if $\mu$ has a density\footnote{If $\mu$ has no density, \eqref{recons} has to be interpreted in the weak sense \eqref{swesk}.}, this density is computed as the discontinuity on the real axis of the Stieltjes transform:
\beq
\label{recons}\mu(x) = \frac{W_{\mu}(x - {\rm i}0) - W_{\mu}(x + {\rm i}0)}{2{\rm i}\pi}\,\dd x\,.
\eeq
Note that there is a unique function $W(z)$ which is holomorphic in $\mathbb{C}\setminus\mathbb{R}$, has a given discontinuity on $\mathbb{R}$, and behaves likes $1/z$ when $|z| \rightarrow \infty$. Indeed, if $\tilde{W}$ was another such function, then $\tilde{W} - W$ would have no discontinuity on $\mathbb{R}$, hence would be holomorphic in $\mathbb{C}$. The growth condition implies that it decays at infinity, and by Liouville theorem, this implies that $\tilde{W} - W = 0$. 

\vspace{0.2cm}

Let us see how it works on a few examples.

\vspace{0.2cm}

\noindent $\bullet$ The Stieltjes transform of a Dirac mass located at $x_0$ is:
$$
W(z) = \frac{1}{z - x_0}\,.
$$
More generally, a simple pole at $z = x_0 \in \mathbb{R}$ with residue $r$ in $W_{\mu}(z)$ indicated that $\mu$ has a contribution from a Dirac mass $r$ located at $x_0$.

\noindent $\bullet$  For the semi-circle law \eqref{SClaw}, we could use the definition \eqref{Stieldef} and compute the integral with the change of variable $x = \sigma(\zeta + 1/\zeta)$ and complex analysis tricks. But there is a better way, relying on \eqref{recons}. Indeed, we are looking for a holomorphic function behaving like $1/z$ when $|z| \rightarrow \infty$, which has a discontinuity on $[-2\sigma,2\sigma]$ such that:
$$
\forall x \in [-2\sigma,2\sigma],\qquad W(x + {\rm i}0) - W(x - {\rm i}0) = -\frac{\sqrt{x^2 - 4\sigma^2}}{\sigma^2}\,.
$$
But we know that the squareroot takes a minus sign when one crosses the locus $[-2\sigma,2\sigma]$ where the quantity inside is negative, so its discontinuity is twice the squareroot. Therefore, the function $-\frac{1}{2\sigma^2}\,\sqrt{z^2 - 4\sigma^2}$ has the discontinuity we look for. It cannot be the final answer for $W(z)$, because of the condition $W(z) \sim 1/z$ when $|z| \rightarrow \infty$. But this can be achieved by adding a polynomial: it does not affect the holomorphicity and discontinuity, but can compensate the growth of the squareroot at infinity. One can check that:
\beq
\label{WSClaw} W_{{\rm sc}}(z) = \frac{z - \sqrt{z^2 - 4\sigma^2}}{2\sigma^2}
\eeq
has all the required properties, provided we choose the determination of the squareroot such that $\sqrt{z^2 - 4\sigma^2} \sim z$ when $|z| \rightarrow \infty$. By uniqueness, \eqref{WSClaw} must be the Stieltjes transform of $\mu_{{\rm sc}}$.

\noindent $\bullet$ Inspired by these two examples, the reader can show that the Stieltjes transform of the Mar\v{c}enko-Pastur law is:
$$
W_{{\rm MP}}(z) = \frac{(1 - \gamma)\sigma^2 + \gamma z - \gamma\sqrt{(z - a_+(\gamma))(z - a_-(\gamma))}}{2\sigma^2 z}\,,
$$
where the determination of the squareroot is fixed by requiring that:
$$
\sqrt{(z - a_+(\gamma))(z - a_-(\gamma))} \sim z
$$
when $|z| \rightarrow \infty$.

\subsection{$\mathcal{R}$-transform}

A closely related tool is the $\mathcal{R}$-transform. To simplify, we consider only measures $\mu$ for which the moments $m_k = \mu[x^k]$ exist for all $k \geq 0$. Let us consider the formal Laurent series:
\beq
\label{Wcurl}\mathcal{W}_{\mu}(z) = \frac{1}{z} + \sum_{k \geq 1} \frac{m_k}{z^{k + 1}}\,.
\eeq
We shall use curly letters to distinguish the formal series from the holomorphic function $W_{\mu}(z)$. There exists a unique formal series:
\beq
\label{Rcurl}\mathcal{R}_{\mu}(w) = \frac{1}{w} + \sum_{\ell \geq 1} \kappa_{\ell}\,w^{\ell - 1}
\eeq
such that:
\beq
\label{RWcurl} \mathcal{R}_{\mu}(\mathcal{W}_{\mu}(z)) = z\,.
\eeq
In other words, $\mathcal{R}_{\mu}$ is the functional inverse -- at the level of formal series -- of $\mathcal{W}_{\mu}$. So, we also have equivalently $\mathcal{W}_{\mu}(\mathcal{R}_{\mu}(w)) = 0$. If we declare that $m_k$ has degree $k$, the $\kappa_{\ell}$ are homogeneous polynomials of degree $\ell$ in the $(m_k)_{k \geq 1}$. One can compute them recursively by replacing \eqref{Wcurl}-\eqref{Rcurl} in \eqref{RWcurl}:
\bea
\kappa_{1} & = & m_1\,, \nonumber \\
\kappa_{2} & = & m_2 - m_1^2\,, \nonumber \\
\kappa_{3} & = & m_3 - 3m_1m_2 + 2m_1^3\,, \nonumber \\
\kappa_{4} & = & m_4 - 4m_1m_3 - 2m_2^2 + 10m_2m_1^2 - 5m_1^4\,,\ldots \nonumber
\eea
The $\kappa_{\ell}$ are called \textbf{free cumulants}. They should not be confused with the better known \textbf{cumulants} $(c_{\ell})_{\ell \geq 1}$, defined by:
$$
\ln\Big(1 + \sum_{k \geq 1} \frac{m_k\,t^{k}}{k!}\Big) = \sum_{\ell \geq 1} \frac{c_{\ell}\,t^{\ell}}{\ell!},\qquad t \rightarrow 0
$$
We see on the first few values:
\bea
c_{1} & = & m_1\,, \nonumber \\
c_{2} & = & m_2 - m_1^2\,, \nonumber \\
c_{3}  & = & m_3 - 3m_1m_2 + 2m_1^3\,, \nonumber \\
c_{4} & = & m_4 - 4m_1m_3 - \textcolor{red}{3}m_2^2 + \textcolor{red}{12}m_2m_1^2 - \textcolor{red}{6}m_1^4\,,\ldots \nonumber
\eea
that $c_{2} = \kappa_{2}$ and $c_{3} = \kappa_{3}$, but this is accidental and in general the cumulants and free cumulants differ for $\ell \geq 4$.

\subsection{Asymptotic freeness}

In general, if $A$ and $B$ are two hermitian matrices, the knowledge of the spectrum of $A$ and $B$ is not enough to determine the spectrum of $A + B$ or $A \cdot B$. Indeed, when $A$ and $B$ do not commute, they cannot be diagonalized in the same basis.

It turns out that for large random matrices "in general position", knowing the spectrum of $A$ and $B$ is enough to reconstruct the spectrum of $A + B$, and the answer is elegantly expressed in terms of the $\mathcal{R}$-transform; the theory is mainly due to Voiculescu around 1991 \cite{Voiculescu}, in the more general context of C$^*$ algebras.  Explaining why this is true would bring us too far, but we aim at presenting the recipe, and illustrating some of its consequences.

We start by introducing several notions, first in a non-random context.
\begin{definition}
If $(M_n)_{n}$ is a sequence of hermitian matrices of size $n$, we say that it has a limit distribution if there exists a probability measure $\mu$ with compact support such that $L^{(M_n)}$ converges to $\mu$ for the vague topology.
\end{definition}

\begin{definition}
\label{AFdef}
Let $(A_n)_{n}$ and $(B_n)_{n}$ two sequences of hermitian matrices of size $n$, admitting as limit distributions respectively $\mu_{A}$ and $\mu_{B}$. We say that $(A_n)_{n}$ and $(B_{n})_{n}$ are \textbf{asymptotically free} if for any positive integers $r,m_1,m_1',\ldots,m_r,m_r'$, we have:
\beq
\label{condor}\lim_{n \rightarrow \infty} n^{-1}\,\mathrm{Tr}\Big\{\prod_{i = 1}^r (A_n^{m_i} - \mu_{A}[x^{m_i}]\cdot I_{n})(B_n^{m_i'} - \mu_{B}[x^{m_i'}]\cdot I_{n})\Big\} = 0\,,
\eeq
where $I_n$ is the identity matrix of size $n$, and the factors in the product are written from the left to the right with increasing $i$.
\end{definition}
If we expand \eqref{condor} and use it recursively, it implies that for asymptotically free matrices, the large $n$ limit of the trace of arbitrary products of $A_n$ and $B_n$ can be computed solely in terms of the moments of $\mu_{A}$ and $\mu_{B}$. In particular, the large $n$ limit of $n^{-1}\mathrm{Tr}\,(A_n + B_n)^{m}$ or $n^{-1}\mathrm{Tr}\,(A_n\cdot B_n)^m$ can be computed solely in terms of $\mu_{A}$ and $\mu_{B}$. Since measures with compact support are determined by their moments, we therefore understand that $\mu_{A}$ and $\mu_{B}$ should determine $\mu_{A + B}$ and $\mu_{A\cdot B}$. Finding the explicit formulas requires some combinatorial work. Focusing on the spectrum of the sum, the result is:
\begin{theorem}
\label{sumth}
If $(A_n)_{n}$ and $(B_n)_{n}$ are asymptotically free and have limit distributions $\mu_{A}$ and $\mu_{B}$, then $(A_n + B_n)_{n}$ has a limit distribution $\mu_{A + B}$, characterized by:
\beq
\label{RABeq} \mathcal{R}_{\mu_{A + B}}(w) = \mathcal{R}_{\mu_{A}}(w) + \mathcal{R}_{\mu_{B}}(w) - \frac{1}{w}\,.
\eeq
\end{theorem}
The last term $-\frac{1}{w}$ is there to ensure that the right-side is of the form $1/w + O(1)$ when $w \rightarrow 0$.

The relevance of this result in random matrix theory is illustrated by the following theorem of Voiculescu:
\begin{theorem}
\label{thV}Let $(A_n)_{n}$ and $(B_n)_{n}$ be two sequences of hermitian random matrices of size $n$. Assume that, for any $n$, $A_n$ is independent of $B_n$, and for any unitary matrix $\Omega_n$, $\Omega_n^{-1} A_n \Omega_n$ is distributed like $A_n$. Then, $(A_n)_{n}$ and $(B_n)_{n}$ are almost surely asymptotically free.
\end{theorem}
In particular, if $L^{(A_n)}$ (resp. $L^{(B_n)}$) converges almost surely to a deterministic $\mu_{A}$ (resp $\mu_{B}$) for the vague topology, using Stieltjes continuity theorem, one deduces that $L^{(A_n + B_n)}$ converges almost surely to a deterministic $\mu_{A + B}$ characterized by \eqref{RABeq}. To compute it, one has to compute the Stieltjes transforms $\mathcal{W}_{\mu_{A}}$ and $\mathcal{W}_{\mu_{B}}$, then compute their functional inverses $\mathcal{R}_{\mu_{A}}$ and $\mathcal{R}_{\mu_{B}}$, use \eqref{RABeq}, compute again the functional inverse $\mathcal{W}_{\mu_{A + B}}$, and finally reconstruct the measure $\mu_{A + B}$ from \eqref{recons}.

\subsection{The semi-circle law as a non-commutative CLT}

From Voiculescu's result, one can understand that the semi-circle law is an analog, in the non-commutative world, of the Gaussian distribution arising when summing independent, identically distributed (i.i.d) real-valued random variables.

Let $(A_n^{(j)})_{1 \leq j \leq N}$ be i.i.d, centered random matrices, whose distribution is invariant under conjugation by a unitary matrix. We assume that the empirical measure of $A_n^{(1)}$ converges almost surely to $\mu_{A}$ for the vague topology. It follows from a slight generalization of Voiculescu's theorem that the family $((A_n^{(j)})_{1 \leq j \leq N})_{n}$ is asymptotically free -- this is defined like in Definition~\ref{AFdef}, except that one uses arbitrary sequences of letters $A^{(j_1)}\cdots A^{(j_s)}$ with $j_i \neq j_{i + 1}$ instead of arbitrary sequences of letters $ABABAB\cdots$. Let us consider:
$$
S_n^{(N)} = \frac{1}{\sqrt{N}} \sum_{j = 1}^N A_n^{(j)}\,.
$$
Theorem~\ref{sumth} has an obvious generalization to this case: for any $N \geq 1$, $(S_n^{(N)})_{n}$ has a limit distribution $\mu_{S^{(N)}}$ when $n \rightarrow \infty$, which is characterized by:
$$
\mathcal{R}_{\mu_{S^{(N)}}} = N\mathcal{R}_{\mu_{A/\sqrt{N}}}(w) - \frac{N - 1}{w}\,.
$$
Playing with the functional equation \eqref{RWcurl}, one easily finds what is the effect of a rescaling on the $\mathcal{R}$-transform:
$$
\mathcal{R}_{\mu_{A}/\sqrt{N}}(w) = N^{-1/2} \mathcal{R}_{\mu_{A}}(N^{-1/2}w)\,.
$$
Since $A_n^{(1)}$ is centered, the first moment of $\mu_{A}$ vanishes. Denoting $\sigma^2_{A}$ the variance of $\mu_{A}$, we can write:
$$
\mathcal{R}_{\mu_{A}}(w) = \frac{1}{w} + \sigma^2\,w + \sum_{\ell \geq 2} \kappa_{\ell + 1}\,w^{\ell}\,,
$$
and therefore:
\beq
\label{Rtru}\mathcal{R}_{\mu_{S^{(N)}}}(w) = \frac{1}{w} + \sigma^2 w + \sum_{\ell \geq 2} N^{(1 - \ell)/2}\,\kappa_{\ell + 1} w^{\ell} \mathop{\longrightarrow}_{N \rightarrow \infty} \frac{1}{w} + \sigma^2 w
\eeq
The functional inverse of $\mathcal{R}_{\infty}(w) = \frac{1}{w} + \sigma^2\,w$ can be readily computed as it is solution of a quadratic equation:
$$
\mathcal{R}_{\infty}(\mathcal{W}_{\infty}(z)) = z \qquad \Longleftrightarrow \qquad \mathcal{W}_{\infty}(z) = \frac{z - \sqrt{z^2 - 4\sigma^2}}{2\sigma^2z}\,.
$$
Note that the determination of the squareroot is fixed by requiring that the formal series $\mathcal{W}_{\infty}(z)$ starts with $1/z + O(1/z)$. We recognize the Stieltjes transform \eqref{WSClaw} of the semi-circle law $\mu_{{\rm sc}}$ with variance $\sigma^2$. Using Stieltjes continuity theorem, one can deduce that $\mu_{S^{(N)}}$ converges for the vague topology to $\mu_{{\rm sc}}$ when $N \rightarrow \infty$. It is remarkable that the limit distribution for $S^{(N)}_{n}$ when $n,N \rightarrow \infty$ does not depend on the details of the summands $A_n^{(j)}$.

Actually, the mechanism of the proof is similar to that of the central limit theorem, provided one replaces the notion of Fourier transform (which is multiplicative for sum of independent real-valued random variables) with the notion of $\mathcal{R}$-transform (which is additive for the sum asymptotically free random matrices). In both cases, the universality of the result -- as well as the occurrence of the Gaussian distribution/the semi-circle law -- comes from the fact that, when the number of summands $N$ goes to infinity, only the second order survives in the formula characterizing the distribution.

\subsection{Perturbation by a finite rank matrix}

We show\footnote{The example we present is inspired by Bouchaud.} how simple computations with the $\mathcal{R}$-transform give insight into the effect of a finite rank perturbation on the spectrum of a GUE matrix. This gives a good qualitative idea of the effect of perturbations on more general random matrices. We will state in Section~\ref{PertW} a complete theorem for Wishart matrices.

So, let $A_n$ be a GUE matrix of size $n$ with variance $\sigma^2$, and consider:
$$
S_n = A_n + B_n,\qquad B_n = {\rm diag}(\underbrace{\Lambda,\ldots,\Lambda}_{m\,\,{\rm times}},\underbrace{0,\ldots,0}_{n - m\,\,{\rm times}})
$$
for $\Lambda > 0$.
We set:
$$
\epsilon = \frac{m}{n}
$$
and would like the study the limit where $n \rightarrow \infty$, and then $\epsilon$ is small. As we have seen, the distribution of $A_n$ is invariant under conjugation by a unitary matrix, and it has the semi-circle law as limit distribution. $B_n$ is deterministic, therefore independent of $A_n$, and it admits a limit distribution given by:
\beq
\label{Bmes}\mu_{B} = (1 - \epsilon)\delta_{0} + \epsilon\delta_{\Lambda}\,.
\eeq
This falls in framework of Voiculescu's theorem, so $S_n$ has a limit distribution $\mu_{S}$. To compute it, we first write down the Stieltjes transform:
$$
W_{\mu_{B}}(z) = \frac{1 - \epsilon}{z} + \frac{\epsilon}{z - \Lambda}\,,
$$
and solving for the functional inverse:
$$
\mathcal{R}_{\mu_{B}}(w) = \frac{1}{2}\bigg[\frac{1}{w} + \Lambda + \sqrt{\Big(\frac{1}{w} - \Lambda)^2 + \frac{4\epsilon\Lambda}{w}}\bigg]\,.
$$
Therefore, we add to it the $\mathcal{R}$-transform \eqref{Rtru} of the semi-circle law minus $1/w$, and we can expand when $\epsilon \rightarrow 0$:
\bea
\mathcal{R}_{\mu_{S}}(w) & = & \sigma^2\,w + \frac{1}{2}\bigg[\frac{1}{w} + \Lambda + \sqrt{\Big(\frac{1}{w} - \Lambda)^2 + \frac{4\epsilon\Lambda}{w}}\bigg] \nonumber \\
& = & \frac{1}{w} + \sigma^2\,w + \frac{\epsilon\Lambda}{1 - \Lambda w} + O(\epsilon^2)\,.
\eea
The Stieltjes transform of $\mu_{S}$ will satisfy:
\beq
\label{eqtosolve}z = \frac{1}{W_{\mu_{S}}(z)} + \sigma^2\,W_{\mu_{S}}(z) +  \frac{\epsilon \Lambda}{1 - \Lambda W_{\mu_{S}}(z)} + O(\epsilon^2)\,.
\eeq
At leading order in $\epsilon$, $\mu_{S}$ the semi-circle law. Let us have a look at the first subleading correction. Qualitatively, two situations can occur.

\noindent $\bullet$ If $W_{{\rm sc}}(z) = 1/\Lambda$ admits a solution $z = z_{\Lambda}$ on the real axis outside of the support $K_{\sigma} = [-2\sigma,2\sigma]$ of $\mu_{S}$, the $O(\epsilon)$ correction to $W_{\mu_{S}}$ has a singularity outside $K_{\sigma}$, which is the sign that $\mu_{S}$ has some mass outside $K_{\sigma}$. If such a real-valued $z_{\Lambda}$ exists, we must have:
$$
\frac{1}{\Lambda} = \frac{z_{\Lambda} - \sqrt{z_{\Lambda}^2 - 4\sigma^2}}{2\sigma^2} \leq \frac{z_{\Lambda} - (z_{\Lambda} - 2\sigma)}{2\sigma^2} \leq \frac{1}{\sigma}\,.
$$
Conversely, if the condition $\Lambda >\sigma$ is met, then there exists a unique such $z_{\Lambda}$, given by:
$$
z_{\Lambda} = \Lambda + \frac{\sigma^2}{\Lambda}\,.
$$
One can then show solving \eqref{eqtosolve} perturbatively that $W_{\mu_{S}}(z)$ has a simple pole at $z = z_{\Lambda} + o(1)$, with residue $\epsilon + o(\epsilon)$. This means that $\mu_{S}$ has a Dirac mass $\epsilon$ at $z_{\Lambda}$. In other words, if $\Lambda$ is above the threshold $\sigma$, a fraction $\epsilon$ of eigenvalues -- i.e. $m = {\rm rank}(B_n)$ eigenvalues -- detach from the support. Even for $\epsilon$ arbitrarily small but non-zero, the maximum eigenvalue is now located at $z_{\Lambda} > 2\sigma$ instead of $2\sigma$ for a GUE matrix.

\noindent $\bullet$ If $\Lambda \leq \sigma$, the singularities of $W_{\mu_{S}}(z)$ remain on $K_{\sigma}$, and therefore the density of $\mu_{S}$ is a small perturbation of the semi-circle, not affecting the position of the maximum eigenvalue.

One should note that the value of the threshold $\Lambda_* = \sigma$ is located in the bulk of the support. We will justify in Section~\ref{FUm} the loose statement that:
\begin{center}
\emph{"eigenvalues of random matrices repel each other"}
\end{center}
This allows an interpretation of the above phenomenon. If we try to add to a random matrix a deterministic matrix with $m$ eigenvalues $\Lambda$, they will undergo repulsion of the eigenvalues that were distributed according to the distribution of $A$ (here, the semi-circle). If the $m$ $\Lambda$'s feel too many eigenvalues of $A$ to their left -- here it happens precisely when $\Lambda > \sigma$ -- they will be kicked out from the support, to a location $z_{\Lambda}$ further to the right of the support. If $\Lambda < \sigma$, the $\Lambda$'s feel the repulsion of enough eigenvalues to their right and to their left to allow for a balance, and thus we just see a small deformation of the semi-circle law, keeping the same support in first approximation.

\newpage

\section{Wishart matrices with perturbed covariance}
\label{PertW}

The same phenomenon was analyzed for complex Wishart matrices by Baik, Ben Arous and P\'ech\'e \cite{BBPphaset}, and is now called the \textbf{BBP phase transition}. The result also holds for real Wishart matrices \cite{VirBloe}. We consider a Wishart matrix $M$ of size $p$, with $n$ degrees of freedom, and covariance $K = {\rm diag}(\Lambda^2,\sigma^2,\ldots,\sigma^2)$. This is a perturbation of the null model with covariance ${\rm diag}(\sigma^2,\ldots,\sigma^2)$.

\begin{theorem}
Assume $n,p \rightarrow \infty$ while $n/p$ converges to $\gamma$, and define:
$$
\Lambda_* = \sigma(1 + \gamma^{-1/2})\,.
$$
\begin{itemize}
\item[$\bullet$] If $\Lambda \in (0,\Lambda_*)$, Theorem~\ref{Fluctumax} continues to hold: $\lambda_{1}^{(M)}$ converges almost surely to $a_+(\gamma)$, and the fluctuations at scale $p^{-2/3}$ follow the Tracy-Widom law.
\item[$\bullet$] If $\Lambda \in (\Lambda_*,+\infty)$, we have almost sure convergence of the maximum:
$$
\lambda_{1}^{(M)} \longrightarrow z_{\Lambda}:= \sigma\Lambda\Big(1 + \frac{\sigma}{\gamma(\Lambda - \sigma)}\Big)\,,
$$
and the random variable
$$
\frac{p^{1/2}}{\sigma\Lambda}\Big(\frac{1}{\gamma} - \frac{\sigma^2}{\gamma^2(\Lambda - \sigma)^2}\Big)^{1/2}\big\{\lambda_{1}^{(M)} - z_{\Lambda}\big\}
$$
describing fluctuations at scale $p^{-1/2}$, converges in law to a Gaussian with variance $1$.
\end{itemize}
\end{theorem}
When $\Lambda$ approaches $\Lambda_*$ at a rate depending on $p$, the maximum eigenvalue converges to $a_+(\gamma)$, but its fluctuations follow a new distribution, that interpolates between Tracy-Widom and Gaussian laws.

For application in statistics, $\Lambda$ can be thought as a trend in empirical data. One may wonder if the trend can be identified from a PCA analysis. The theorem shows that the answer is positive only if the trend is strong enough -- i.e. $\Lambda > \Lambda_*$. As for perturbation of the GUE, the threshold $\Lambda_*^2$ lies inside the support of the Mar\v{c}enko-Pastur law.

Although more interesting for statistics, the case of real Wishart matrices was only tackled in 2011 by Bloemendal and Vir\'ag\footnote{Actually, their method relate the distributions for the fluctuations of the maximum of perturbed GOE or GUE to the probability of explosion of the solution of second order stochastic differential equation. In the unperturbed case, they also obtained characterizations of the same nature for the Tracy-Widom laws. This is a beautiful result fitting in the topic of the summer school, however at a more advanced level compared to the background provided at the school.}, with similar conclusions. The reason is that, in the complex case, we will see in Section~\ref{Exact} that algebraic miracles greatly facilitates the computations, which boil down to analyzing the asymptotic behavior of a sequence of orthogonal polynomials. This can be done with the so-called Riemann-Hilbert steepest descent analysis, developed by Deift, Zhou and coauthors in the 90s -- for an introduction, see \cite{Defcours} -- and this is the route taken by BBP.

\newpage

\section{From matrix entries to eigenvalues}
\label{Eig}

\subsection{Lebesgue measure and diagonalization}
\label{FUm}
We would like to compute the joint distribution of eigenvalues of a symmetric or hermitian random matrix. For this purpose, we basically need to perform a change of variables in integrals of the form $\int \dd M\,f(M)$, hence to compute the determinant of the Jacobian of this change of variable. Although some details have to be taken care of before arriving to that point, the core of the computation is easy and concentrated in \eqref{diffM} and the evaluation of the determinant.

First consider the case of symmetric matrices. Let $\mathcal{O}_{n}$ be the set of orthogonal $n \times n$ matrices, i.e. satisfying $\Omega^{T}\Omega = I_{n}$. Since any symmetric matrix can be diagonalized by an orthogonal matrix, the $\mathcal{C}^{\infty}$ map:
\beq
\label{ChangeofS}M\,:\,\begin{array}{ccc} \mathcal{O}_{n} \times \mathbb{R}^n &  \longrightarrow & \mathcal{S}_{n} \\
(\Omega,\lambda_1,\ldots,\lambda_n) & \longmapsto & \Omega\,{\rm diag}(\lambda_1,\ldots,\lambda_n) \Omega^{-1}\end{array}
\eeq
is surjective. However, the map is not injective, so we cannot take \eqref{ChangeofS} as an admissible change of variable.
Indeed, if:
$$
M = \Omega {\rm diag}(\lambda_1,\ldots,\lambda_n)\Omega^{-1} = \tilde{\Omega}{\rm diag}(\tilde{\lambda}_1,\ldots,\tilde{\lambda}_n)\tilde{\Omega}^{-1},
$$
then there exists a permutation $\sigma \in \mathfrak{S}_{n}$ and an orthogonal matrix $D$ that leaves stable the eigenspaces of $M$ such that:
\beq
\label{Dun}\tilde{\Omega} = \Omega D,\qquad \tilde{\lambda}_i = \lambda_{\sigma(i)}\,.
\eeq
To solve this issue, we first restrict to the subset $(\mathcal{S}_{n})_{\Delta}$ consisting of symmetric matrices with pairwise distinct eigenvalues. This is harmless since $(\mathcal{S}_{n})_{\Delta}$ is an open dense subset of $\mathcal{S}_{n}$, hence its complement has Lebesgue measure $0$. Then, two decompositions are related by \eqref{Dun} with $D$ being a diagonal orthogonal matrix, and this forces the diagonal entries to be $\pm 1$. So, let us mod out the left-hand side of \eqref{ChangeofS} by $\{\pm 1\}^n$. Then, we can kill the freedom of permuting the $\lambda_i$'s by requiring that $\lambda_i$ decreases with $i$. Denoting:
$$
(\mathbb{R}_{n})_{\Delta} = \big\{(\lambda_1,\ldots,\lambda_n) \in \mathbb{R}^n,\qquad \lambda_1 > \lambda_2 > \ldots >\lambda_n\big\}\,,
$$
we finally obtain an invertible map:
\beq
\label{ChangeofS2}M\,:\,\begin{array}{ccc} \big(\mathcal{O}_{n}/\{\pm 1\}^n\big) \times (\mathbb{R}^n)_{\Delta} &  \longrightarrow & (\mathcal{S}_{n})_{\Delta} \\
(\Omega,\lambda_1,\ldots,\lambda_n) & \longmapsto & \Omega\,{\rm diag}(\lambda_1,\ldots,\lambda_n) \Omega^{-1} \end{array}
\eeq
and one can show that it is a $\mathcal{C}^{\infty}$ diffeomorphism -- i.e. an admissible change of variable.

To be more explicit, we have to choose coordinates on $\mathcal{O}_{n}$. In the vicinity of $I_n \in \mathcal{O}_{n}$, we can choose as coordinates the entries $(\omega_{ij})_{1 \leq i < j \leq n}$ of an antisymmetric matrix $\omega$, which parametrizes an orthogonal matrix by the formula $\Omega = \exp(\omega)$. And in $\mathcal{S}_{n}$, we remind that we had chosen as coordinates the entries $(M_{ij})_{1 \leq i \leq j \leq n}$. Then, we know that:
$$
\dd M = 2^{-n}\,\prod_{1 \leq i < j \leq n} \!\!\!\! \dd \omega_{ij}\,\prod_{i = 1}^n \dd \lambda_i\,\mathcal{J}(\lambda,\omega)\,,
$$
where the $2^{-n}$ comes from the quotient by $\{\pm 1\}^{n}$, and it remains to compute the Jacobian determinant:
\begin{figure}[h!]
\begin{center}
\includegraphics[width=0.40\textwidth]{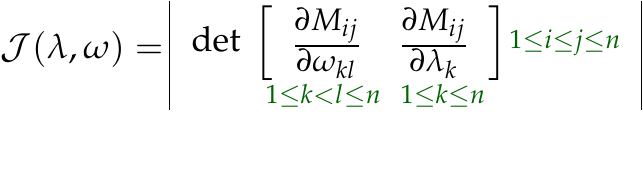}
\end{center}
\end{figure}

First, we remind that the Lebesgue measure is invariant under conjugation of $M$ by an orthogonal matrix.  We can thus evaluate the derivatives at $\omega = 0$ (i.e. $\Omega = I_n$) and find:
\beq
\label{diffM} \dd M_{ij} = [\dd \omega,\Lambda]_{ij} + \dd \Lambda_i \delta_{ij} = \dd \omega_{ij}(\lambda_i - \lambda_j) + \dd\lambda_i \delta_{ij}
\eeq
Therefore, the matrix in the Jacobian is diagonal: in the first block $1 \leq i < j \leq n$ and $1 \leq k < l \leq n$, the diagonal elements $(i,j) = (k,l)$ are $(\lambda_i - \lambda_j)$, and in the second block, the diagonal elements are just $1$. Therefore:
$$
\mathcal{J}(\lambda,0) = \prod_{1 \leq i < j \leq n} |\lambda_j - \lambda_i|
$$

We can repeat all steps for hermitian matrices. $\mathcal{O}_{n}$ should be replaced with the set $\mathcal{U}_{n}$ of unitary matrices, i.e. satisfying $(\Omega^{T})^*\Omega = I_n$. The map \eqref{ChangeofS} now sends $\mathcal{U}_{n} \times \mathbb{R}^n$ to $\mathcal{H}_{n}$. It is not surjective, but if we restrict to the set $(\mathcal{H}_{n})_{\Delta}$ of hermitian matrices with pairwise distinct eigenvalues, the only freedom is to have \eqref{Dun} with $D$ a diagonal matrix whose entries are complex numbers of unit norm ; we denote $\mathcal{U}_{1}^n$ the group of such matrices. Then, we obtain an admissible change of variable:
\beq
\label{ChangeofH2} \big(\mathcal{U}_{n}/\mathcal{U}_{1}^n\big) \times (\mathbb{R}^n)_{\Delta} \simeq (\mathcal{H}_{n})_{\Delta}\,.
\eeq
As coordinates on $\mathcal{U}_{n}$ near $I_n$, we can take the real and imaginary parts of the entries $(\omega_{ij})_{1 \leq i < j \leq n}$ of a matrix $\omega$ such that\footnote{Such a matrix is called "antihermitian".}  $(\omega^{T})^* = -\omega$, parametrizing a unitary matrix by the formula $\Omega = \exp(\omega)$. The formula \eqref{diffM} for the differential does not change but we have now twice many coordinates: the Jacobian matrix is still diagonal, and the diagonal entries corresponding to derivative with respect to ${\rm Re}\,\omega_{ij}$ and to ${\rm Im}\,\omega_{ij}$ both evaluate to $(\lambda_i - \lambda_j)$. Thus, the Jacobian determinant reads:
$$
\mathcal{J}(\lambda,0) = \prod_{1 \leq i < j \leq n} |\lambda_j - \lambda_i|^{2}\,.
$$

There is a last step about which we will be brief: this result -- valid at $\Omega = I_n$ -- has to be transported to any point of $\mathcal{S}_{n}$ (or $\mathcal{H}_{n})$ by conjugating with an $\mathcal{O}_{n}$ (resp $\mathcal{U}_{n}$) matrix. Of course, this does not affect the eigenvalue dependence of the Jacobian factor. The result makes appear the Haar measure on $\mathcal{O}_{n}$ (resp. $\mathcal{U}_{n}$): this is the unique probability measure which is invariant under left and right multiplication by an orthogonal (resp. unitary) matrix. We denote $\dd\nu(\Omega)$ the measure induced by the Haar measure on the quotient $\mathcal{O}_{n}/\{\pm 1\}^n$ (resp. $\mathcal{U}_{n}/\mathcal{U}_{1}^n$).

\begin{theorem}
\label{WeylL} Under the change of variable \eqref{ChangeofS2} or \eqref{ChangeofH2}, we have:
$$
\dd M = c_{\beta,n}\,\dd\nu(\Omega)\,\prod_{i = 1}^n \dd\lambda_i \prod_{1 \leq i < j \leq n} |\lambda_j - \lambda_i|^{\beta}
$$
for some (explicitly computable) constant $c_{\beta,n} > 0$.
\end{theorem}

\subsection{Repulsion of eigenvalues}

As a consequence, if $M$ is a random symmetric (resp. hermitian) matrix whose p.d.f. of entries is $\dd M\,F(M)$, and $f$ is invariant under conjugation by an orthogonal (resp. unitary) matrix, then $F(M)$ is actually a function $f(\lambda_1,\ldots,\lambda_n)$ of the eigenvalues only, and the joint p.d.f of the eigenvalues of $M$ is proportional to:
\beq
\label{jointpdf} Z_{n,\beta}^{-1} \prod_{1 \leq i < j \leq n} \big|\Delta(\lambda_1,\ldots,\lambda_n)\big|^{\beta}\,f(\lambda_1,\ldots,\lambda_n)\,,
\eeq
with:
\beq
\label{Delt}\Delta(\lambda_1,\ldots,\lambda_n) = \prod_{1 \leq i < j \leq n} (\lambda_j - \lambda_i)\,,
\eeq
and the constant $Z_{n,\beta}$ is such that the integral of \eqref{jointpdf} against the Lebesgue measure over $\mathbb{R}^n$ evaluates to $1$. Because of the factor $|\Delta(\lambda_1,\ldots,\lambda_n)|^{\beta}$ the probability that two eigenvalues are close to each other is small: the eigenvalues of a random matrix usually repel each other. The intensity of the repulsion is measured by the parameter $\beta$, which is fixed by the type of the matrix (symmetric or hermitian).

\begin{lemma}
\eqref{Delt} is the \textbf{Vandermonde determinant}:
$$
\Delta(\lambda_1,\ldots,\lambda_n) = \det\left[\begin{array}{cccc} 1 & 1 & \cdots & 1 \\ \lambda_1 & \lambda_2 & \cdots & \lambda_n \\ \vdots & \vdots &  & \vdots \\ \lambda_1^{n - 1} & \lambda_2^{n - 1} & \cdots & \lambda_n^{n - 1} \end{array}\right]\,.
$$
\end{lemma}
\noindent \textbf{Proof.} Let us denote $D(\lambda_1,\ldots,\lambda_n)$ the determinant in the right-hand side. It is a polynomial function of $\lambda_i$, of degree at most $n - 1$, which admits the $n - 1$ roots $\lambda_i = \lambda_j$ indexed by $j \neq i$. Therefore, we can factor out successively all the monomials that occur in $\Delta$, and find:
\beq
\label{umsa} D(\lambda_1,\ldots,\lambda_n) = c_n\,\Delta(\lambda_1,\ldots,\lambda_n)
\eeq
for some constant $c_n$. We prove by induction that $c_n = 1$. This is obviously true for $n = 1$. If this is true for $(n - 1)$, we expand the determinant of size $n$ with respect to its last column, and find that the coefficient of $\lambda_{n}^{n - 1}$ is $D(\lambda_1,\ldots,\lambda_{n - 1})$. Comparing with \eqref{umsa} and the induction hypothesis, we deduce that $c_n = 1$. \hfill $\Box$

\vspace{0.2cm}

\begin{lemma}
\label{VD2s}
For any sequence $(Q_m)_{m \geq 0}$ of polynomials of degree $m$ with leading coefficient $1$:
$$
\Delta(\lambda_1,\ldots,\lambda_n) = \mathop{{\rm det}}_{1 \leq i,j \leq n} \big[Q_{i - 1}(\lambda_j)\big]\,.
$$
\end{lemma}
\noindent \textbf{Proof.} By adding linear combinations of the $(n - 1)$ first lines to the last line, one can actually replace $\lambda_j^{n - 1}$ in the last line by $Q_{n - 1}(\lambda_j)$ for any polynomial $Q_{n - 1}$ of degree $n - 1$ with leading coefficient $1$. Repeating this procedure successively for the lines $(n - 1)$, $(n - 2)$, etc. establishes the claim. \hfill $\Box$

\subsection{Eigenvalue distribution of Wishart matrices}

The result for Wishart matrices was obtained almost simultaneously in 1939 by \cite{WFisher,WGirschik,WHsu,WRoy}.

\begin{theorem}
If $M$ is a real ($\beta = 1$) or complex ($\beta = 2$) Wishart matrix with covariance $K = {\rm diag}(\sigma^2,\ldots,\sigma^2)$, of size $p$ with $n$ degrees of freedom, the joint p.d.f of its eigenvalues is:
\beq
\label{announcW}Z_{n,\beta}^{-1}\,\prod_{1 \leq i < j \leq n} |\lambda_i - \lambda_j|^{\beta}\,\prod_{i = 1}^n \lambda_i^{\frac{\beta}{2}(n - p) + \frac{\beta - 2}{2}}\,\exp\Big(-\frac{n\beta}{2\sigma^2}\,\lambda_i\Big)
\eeq
for an (explicitly computable) normalization constant $Z_{n,\beta}^{-1}$.
\end{theorem}

\noindent \textbf{Proof.} The proof is a bit more involved than in Section~\ref{FUm}, and was omitted during the lectures. It uses a change of variable in three steps, the last one being already given by Theorem~\ref{WeylL}. We give the details for the case of real Wishart matrices.

\vspace{0.1cm}

\noindent $\bullet$ First, we consider $X$ as a matrix of $p$ vectors in $\mathbb{C}^n$, which we can orthogonalize. This produces in a unique way a matrix $\Omega$ of size $n \times p$, such that:
\beq
\label{Omegrel}\Omega^{T}\Omega = I_{p}\,.
\eeq
and a lower  triangular matrix $L$ of size $p \times p$ with positive diagonal entries, such that:
\beq
\label{XOL} X = \Omega L\,.
\eeq
The Lebesgue measure $\dd X$ is invariant under multiplication to the left by an orthogonal matrix of size $n$, thus it is enough to evaluate the Jacobian at $\Omega$ equals:
$$ 
\Omega^{0} = \left[\begin{array}{c} I_{p,p} \\ 0_{n - p,p} \end{array}\right]\,,
$$
where $0_{m,p}$ is the matrix of size $m \times p$ filled with $0$'s. 

We need to fix local coordinates on the tangent space at $\Omega^0$ of the set $\mathcal{O}_{n,p}$ of matrices $\Omega$ satisfying \eqref{Omegrel}. For example, we can choose the entries $\Omega_{kl}$ with $1 \leq k < l \leq p$, and the $\Omega_{kl}$ with $k \geq p + 1$ and $1 \leq l \leq p$. The remaining $\Omega_{kl}$ with $1 \leq l < k \leq p$ are then determined by \eqref{Omegrel}, and infinitesimally around $\Omega^{0}$ we find for these indices $\Omega_{kl} = -\Omega_{lk}$. The dimension of $\mathcal{O}_{n,p}$ is thus $p(p - 1)/2 + p(n - p)$. For the matrix $L$, we naturally choose as coordinates its non-zero entries $L_{kl}$ indexed by $1 \leq l \leq k \leq p$ -- the space of $L$'s has dimension $p(p + 1)/2$. This is consistent with the dimension of the space of $X$'s:
$$
np = \frac{p(p + 1)}{2} + \frac{p(p - 1)}{2} + p(n - p)\,.
$$
Now, we compute the differential of \eqref{XOL}:
$$
\dd X_{ij} = \delta_{ik}\delta_{jl} \dd L_{kl} + \dd\omega_{kl} \delta_{ik} \delta_{k > l} - \dd \omega_{kl} L_{kl} \delta_{il}\delta_{k < l} + \dd \omega_{kl} L_{lj} \delta_{ik} \delta_{k > p}\,.
$$
A careful look at the indices shows that the Jacobian matrix is of the form:
\begin{figure}[h!]
\begin{center}
\includegraphics[width=0.50\textwidth]{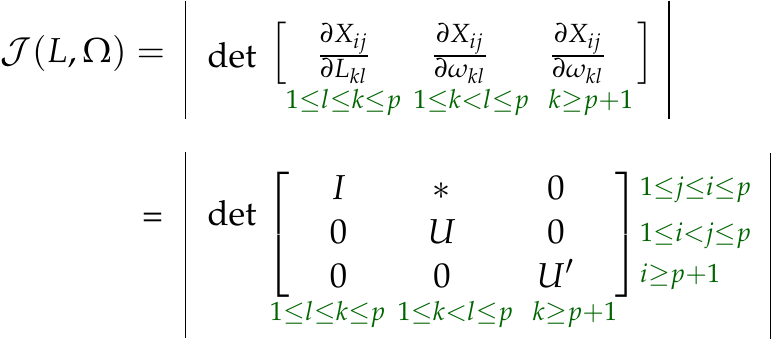} 
\end{center}
\end{figure}

\noindent with $U$ and $U'$ upper triangular matrices with respect to the lexicographic order on the ordered pair $(i,j)$. Besides, the diagonal elements of $U$ and $U'$ at position $(i,j) = (k,l)$ are $L_{jj}$. So, the determinant evaluates to:
$$
\mathcal{J}(L,\Omega) = \prod_{j = 1}^p L_{jj}^{n - p + j - 1}\,,
$$
and we have:
\beq
\label{C1}\dd X = \dd \nu(\Omega)\,\prod_{1 \leq j \leq i \leq p} \dd L_{ij}\,\prod_{j = 1}^p L_{jj}^{n - p + j - 1}\,,
\eeq
where $\dd\nu(\Omega)$ is the measure on $\mathcal{O}_{n,p}$ obtained by transporting the volume element of the $\omega$'s from $\Omega^0$ to any point in $\mathcal{O}_{n,p}$.

\vspace{0.1cm}

\noindent $\bullet$ Next, we change variables from $L$ to $M$:
$$
M = n^{-1}\,X^{T}X = n^{-1}\,L^{T}L\,.
$$
The differential is:
$$
\dd M_{ij} = n^{-1}\Big(\delta_{lj} L_{ki} + \delta_{li} L_{kj}\Big)\,,
$$
and we must compute the Jacobian:
\begin{figure}[h!]
\begin{center}
\includegraphics[width=0.35\textwidth]{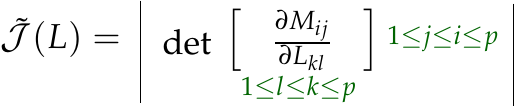}
\end{center}
\end{figure}

If we put on ordered pairs $(i,j)$ the lexicographic order, we observe that the Jacobian matrix is upper triangular, with entries $n^{-1}(\delta_{jj}L_{ii} + \delta_{ij}L_{jj})$ on the diagonal with double index $(i,j)$. Therefore:
\beq
\label{C2}\dd M = \dd L\,\tilde{\mathcal{J}}(L)\,,\qquad \tilde{\mathcal{J}}(L) = n^{-p(p + 1)/2}\,2^{p} \prod_{j = 1}^p L_{jj}^{j}\,.
\eeq

\vspace{0.1cm}

\noindent $\bullet$ Combining \eqref{C1} and \eqref{C2} yields:
$$
\dd X = c_{n,p}\,\dd\nu(\Omega)\,\dd M\,\prod_{j = 1}^p L_{jj}^{n - p - 1}
$$
and we rewrite:
\bea
\prod_{j = 1}^n L_{jj}^{n - p - 1} & = & \det(L)^{n - p - 1} = \det(L^{T}L)^{(n - p - 1)/2} \nonumber \\
& = & n^{p(p + 1)/2}\,\det(M)^{(n - p - 1)/2} = n^{p(p + 1)/2} \prod_{j = 1}^n \lambda_j^{(n - p - 1)/2} \,.\nonumber
\eea
Finally, we use Theorem~\ref{WeylL} to obtain the announced result \eqref{announcW} in the case $\beta = 1$.

\vspace{0.1cm}

\noindent $\bullet$ The case of complex Wishart matrices is treated similarly, with $\mathcal{O}_{n,p}$ being replaced by the set $\mathcal{U}_{n,p}$ of $n \times p$ matrices $\Omega$ such that $(\Omega^{T})^*\Omega = I_{p}$. \hfill $\Box$

\newpage

\section{Exact computations in invariant ensembles}

\subsection{Invariant ensembles}
\label{invens}
The Gaussian ensembles and the Wishart ensembles are special cases of the \textbf{invariant ensembles}. These are symmetric (resp. hermitian) random matrices $M$ of size $n$, whose distribution of entries is of the form:
\beq
\label{DZ}Z_{n,\beta}^{-1}\,\,\dd M\,\exp\Big(-\frac{n\beta}{2}\,\mathrm{Tr}\,V(M)\Big)\,.
\eeq
The function $V$ is assumed to grow fast enough at infinity -- e.g. $V$ is a polynomial with positive leading coefficient -- so that \eqref{DZ} has finite mass on $\mathcal{S}_{n}$ or $\mathcal{H}_{n}$, and we tune $Z_{n,\beta}^{-1}$ so that this mass is $1$. Theorem~\ref{WeylL} implies that the joint p.d.f of the eigenvalues\footnote{Contrarily to the previous sections, in \eqref{VW} the eigenvalues are not assumed to be ordered. When we need to consider the maximum eigenvalue, we shall use the notation $\lambda_{\max}$.} is:
\beq
\label{VW} Z_{n,\beta}^{-1}\,\prod_{1 \leq i < j \leq n} |\lambda_i - \lambda_j|^{\beta}\,\,\prod_{i = 1}^n \exp\Big\{-\frac{n\beta}{2}\,V(\lambda_i)\Big\}\,.
\eeq
The Wishart ensembles -- in which the size is denoted $p$ instead of $n$ -- correspond to the cases:
\beq
\label{WisV}V(x) =  -\frac{x}{\sigma^2} +\Big[\gamma - 1 + \frac{1}{p}\Big(1 - \frac{2}{\beta}\Big)\Big]\ln x,\qquad \gamma = n/p\,,
\eeq
and the Gaussian ensembles to:
$$
V(x) = \frac{x^2}{2\sigma^2}\,.
$$

Note that the distribution \eqref{VW} makes sense for any value of $\beta > 0$. When $\beta$ increases starting from $0$, it provide a model interpolating from independent random variables to strongly correlated (repulsive) random variables, called the \textbf{$\beta$-ensembles}.

Equation~\ref{VW} still contains too much information. We would like to answer questions like: what is the probability that one eigenvalue falls into a given interval ? In other words, we want to compute the marginals of the distribution \eqref{VW}. Surprisingly, for $\beta = 1$ and $\beta = 2$, this can be performed exactly, using tricks mainly discovered by Gaudin and Mehta in the early 60s. We will stick to the case $\beta = 2$, for which the computations are in fact much simpler. And since for the moment we will be occupied with exact computations, it is convenient to use a notation $W(\lambda_i)$ instead of $(n\beta/2)V(\lambda_i)$ in \eqref{VW}.

\subsection{Partition function}
\label{partf}
Prior to any computation, it is useful to evaluate the normalization constant, also called \textbf{partition function}
$$
Z_{n} = \int_{\mathbb{R}^n} \prod_{1 \leq i < j \leq n} \big|\Delta(\lambda_1,\ldots,\lambda_n)\big|^2\,\prod_{i = 1}^n e^{-W(\lambda_i)}\,.
$$
This can be done in terms of the orthogonal polynomials $(P_n)_{n \geq 0}$ for the measure $\dd x\,e^{-W(x)}$ on $\mathbb{R}$. More precisely, consider the scalar product on the space of real-valued polynomials:
\beq
\label{scalar}\langle f,g\rangle = \int_{\mathbb{R}} f(x)\,g(x)\,e^{-W(x)}\,\dd x\,.
\eeq
The orthogonalization of the canonical basis $(x^n)_{n \geq 0}$ for the scalar product \eqref{scalar} determines a unique sequence $(p_n)_{n \geq 0}$ of polynomials with the following properties:
\begin{itemize}
\item[$\bullet$] $P_n$ has degree $n$ and starts with $x^n + \cdots$.
\item[$\bullet$] For any $n,m \geq 0$, $\langle P_n,P_m \rangle = \delta_{nm} h_{n}$ for some constant $h_n > 0$.
\end{itemize}

\begin{theorem}
\label{PZ}
$\boxed{Z_{n} = n! \prod_{m = 0}^{n - 1} h_m}$
\end{theorem}
\noindent \textbf{Proof.} Let $(Q_m)_{m \geq 0}$ be an arbitrary sequence of polynomials of degree $m$ with leading coefficient $1$, use the representation of Lemma~\ref{VD2s} for the Vandermonde determinant, and expand the determinants: 
$$
Z_{n} = \sum_{\sigma,\tau \in \mathfrak{S}_{n}} {\rm sgn}(\sigma){\rm sgn}(\tau) \int_{\mathbb{R}^n}\prod_{i = 1}^n Q_{\sigma(i) - 1}(\lambda_i)\,Q_{\tau(i) - 1}(\lambda_i)\,e^{-W(\lambda_i)}\dd\lambda_i\,. \nonumber
$$
We observe that, in each term, the integral over $\mathbb{R}^n$ factors into $n$ integrals over $\mathbb{R}$. Then, $i$ is a dummy index for the product, and we can also rename it $\tau^{-1}(i)$. Since the signatures satisfy ${\rm sgn}(\sigma){\rm sgn}(\tau) = {\rm sgn}(\sigma\tau^{-1})$, we shall change variables in the sum and set $\tilde{\sigma} = \sigma\tau^{-1}$. The summands only depend on $\tilde{\sigma}$, and it remains a sum over a permutation, which produces a factor of $n!$. So:
\bea
Z_{n} & = & n!\,\sum_{\tilde{\sigma} \in \mathfrak{S}_{n}} {\rm sgn}(\tilde{\sigma}) \prod_{i = 1}^n \Big[ \int_{\mathbb{R}} Q_{\tilde{\sigma}(i) - 1}(x)\,Q_{i - 1}(x)\,e^{-W(x)}\dd x\Big]\nonumber \\
& = & n! \mathop{{\rm det}}_{1 \leq i,j \leq n} \Big[\int_{\mathbb{R}} Q_{i - 1}(x)Q_{j - 1}(x)\,e^{-W(x)}\Big]\,,
\eea
where, in the last line, we have used the multilinearity of the determinant. Now, if we choose $(Q_m)_{m \geq 0}$ to be the orthogonal polynomials for the scalar product \eqref{scalar}, the matrix in the determinant becomes diagonal. This entails the result. \hfill $\Box$

\subsection{Marginals of eigenvalue distributions}

\subsubsection{J\'anossy densities}

If $M$ is a random hermitian matrix, we define the $k$-point \textbf{J\'anossy densities} $\rho_{n}^{(k)}(x_1,\ldots,x_k)$, as the functions such that, for any pairwise disjoint measurable sets $A_1,\ldots,A_k$:
\beq
\label{Janos}\mathbb{P}\Big[\exists i_1,\ldots,i_k,\qquad \lambda_{i_j} \in A_{j}\Big] = \int_{A_1 \times \cdots \times A_k} \rho_{n}^{(k)}(x_1,\ldots,x_k)\,\prod_{i = 1}^k \dd x_i\,.
\eeq

The $\rho_{n}^{(k)}$ can be considered as a probability density -- in particular they are non-negative -- except that their total integral is not $1$. Since the eigenvalues are not ordered in \eqref{Janos}, $\rho_n^{(k)}$ is a symmetric function of $x_1,\ldots,x_k$, and we have:
\beq
\label{normPO}\int_{\mathbb{R}^k} \rho_{n}^{(k)}(x_1,\ldots,x_k) \prod_{i = 1}^k \dd x_i = \frac{n!}{(n - k)!}\,,
\eeq
i.e. the number of ways of choosing $k$ ordered eigenvalues among $n$. The $1$-point J\'anossy density coincides with the average spectral density multiplied by $n$, since
$$
\int_{\mathbb{R}} \rho_{n}^{(1)}(x)\,\dd x = n\,.
$$
Besides, $\rho_{n}^{(n)}$ is nothing but the joint p.d.f of the $n$-eigenvalues, multiplied by $n!$ since \eqref{normPO} gives:
$$
\int_{\mathbb{R}^n} \rho_{n}^{(n)}(x_1,\ldots,x_n) = n!\,.
$$
The $k$-point densities can be found by integrating out $(n - k)$ variables in $\rho_{n}^{(n)}$, again paying attention to the normalization constant:
\beq
\label{rhonk}\rho_{n}^{(k)}(x_1,\ldots,x_k) = \frac{1}{(n - k)!} \int_{\mathbb{R}^{n - k}} \rho_n^{(n)}(x_1,\ldots,x_n)\,\prod_{i = k + 1}^n \dd x_i\,.
\eeq

\subsubsection{In invariant ensembles}

When the random matrix is drawn from an invariant ensemble (Section~\ref{invens}), we have:
\beq
\label{pnn}\rho_{n}^{(n)}(x_1,\ldots,x_n) = \frac{n!}{Z_{n}}\,\Delta(x_1,\ldots,x_n)^2\,\prod_{i = 1}^n e^{-W(x_i)}\,.
\eeq
The J\'anossy densities can be computed in terms of the orthogonal polynomials which already appeared in Section~\ref{partf} to compute $Z_{n}$. Let us introduce the Christoffel-Darboux kernel:
$$
K_n(x,y) = \sum_{k = 0}^{n - 1} \frac{P_k(x)P_k(y)}{h_k}\,.
$$
Using the orthogonality relations, one can easily prove:
\beq
\label{CKD}K_n(x,y) = \frac{P_{n}(x)P_{n - 1}(y) - P_{n - 1}(x)P_{n}(y)}{h_{n - 1}(x - y)}\,,
\eeq
which is more advantageous -- especially from the point of the large $n$ regime -- since it only involves two consecutive orthogonal polynomials.

\begin{theorem}
\beq
\label{DetP} \boxed{\rho_{n}^{(k)}(x_1,\ldots,x_k) = \mathop{{\rm det}}_{1 \leq i,j \leq k} \big[\tilde{K}_{n}(x_i,x_j)\big]}\,,
\eeq
where $\tilde{K}_n(x,y) = K_n(x,y)\,e^{-[W(x) + W(y)]/2}$.
\end{theorem}
\noindent \textbf{Proof.} We first consider $k = n$. With \eqref{pnn} and Lemma~\ref{VD2s} and Theorem~\ref{PZ}, we can write:
$$
\rho_{n}^{(n)}(x_1,\ldots,x_n) = \frac{n!}{n!\,\prod_{m = 0}^{n - 1} h_m} \det_{1 \leq i,j \leq n} \big[P_{j - 1}(x_i)\big]\cdot \det_{1 \leq k,l \leq n} \big[P_{k - 1}(x_l)\big]\,\prod_{i = 1}^n e^{-W(\lambda_i)}\,.
$$
We implicitly used $\det (A^{T}) = \det(A)$ to write the first determinant. We then push a factor $h_m^{1/2}$ in the columns (resp. in the lines) of the first (resp. the second) determinant, and a factor $\exp[-W(\lambda_m)/2]$ in the lines (resp. the columns) of the first (resp. the second) determinant. The result, using $\det (A\cdot B) = (\det A) \cdot (\det B)$, reads:
\bea
\rho_{n}^{(n)}(x_1,\ldots,x_n) & = & \mathop{{\rm det}}_{1 \leq i,j \leq n} \big[h_{j - 1}^{-1/2}P_{j - 1}(x_i)\,e^{-W(x_i)/2}\big] \cdot \mathop{{\rm det}}_{1 \leq k,l \leq n} \big[h_{k - 1}^{-1/2}P_{k - 1}(x_l)\,e^{-W(x_l)/2}\big] \nonumber \\
& = & \mathop{{\rm det}}_{1 \leq i,l \leq n}\Big[\sum_{k = 1}^{n} \frac{P_{k - 1}(x_i)P_{k - 1}(x_l)}{h_{k - 1}}\,e^{-[W(x_i) + W(x_l)]/2}\Big]\,, \nonumber
\eea 
which is the desired result.

Next, we would like to integrate out the last $n - k$ variables in $\rho_{n}^{(n)}$ to find $\rho_{n}^{(k)}$ via \eqref{rhonk}. This is achieved by successive application of the one-step integration lemma:
\begin{lemma}
\beq
\label{RK}\int_{\mathbb{R}} \mathop{{\rm det}}_{1 \leq i,j \leq k}\big[ \tilde{K}_{n}(x_i,x_j)\big]\,\dd x_k = (n - k + 1) \mathop{{\rm det}}_{1 \leq i,j \leq k - 1} \big[\tilde{K}_{n}(x_i,x_j)\big]\,.
\eeq
\end{lemma}

To prove the lemma, we first remark that $\tilde{K}_n(x,y)$ is the kernel of an operator $\widehat{K}_{n}\,:\,L^2(\mathbb{R},\dd x) \longrightarrow L^2(\mathbb{R},\dd x)$, which is the orthogonal projection onto the rank $n$ subspace
$$
V_n = \mathbb{R}_{n - 1}[x]\cdot e^{-W(x)/2}
$$
In particular -- as one can check directly:
\bea
\int_{\mathbb{R}} \tilde{K}_n(x,z)\tilde{K}_{n}(z,y)\,\dd z & = & \tilde{K}_{n}(x,y)\,, \nonumber \\
\int_{\mathbb{R}} \tilde{K}_{n}(z,z)\,\dd z & = & n\,. \nonumber 
\eea
Let us expand the $k \times k$ determinant in the left-hand side of \eqref{RK}:
$$
\int_{\mathbb{R}} \mathop{{\rm det}}_{1 \leq i,j \leq k}\big[\tilde{K}_{n}(x_i,x_j)\big]\,\dd x_k = \sum_{\sigma \in \mathfrak{S}_{k}} {\rm sgn}(\sigma)\,\int_{\mathbb{R}} \Big[\prod_{i = 1}^k \tilde{K}_{n}(x_i,x_{\sigma(i)})\Big]\,\dd x_k\,.
$$
We find two types of terms:
\begin{itemize}
\item[$\bullet$] If $\sigma(k) = k$, we have a factor $$\int_{\mathbb{R}} \tilde{K}_{n}(x_k,x_k)\,\dd x_k = n.$$ The remaining factors is a sum over all permutations $\tilde{\sigma} \in \mathfrak{S}_{k - 1}$, which reconstructs
$$
\mathop{{\rm det}}_{1 \leq i,j \leq k - 1}\big[\tilde{K}_{n}(x_i,x_j)\big]\,.
$$
\item[$\bullet$] If $\sigma(k) \neq k$, we rather have a factor $$\int_{\mathbb{R}} \tilde{K}_{n}(x_{\sigma^{-1}(k)},x_k)\,\tilde{K}_{n}(x_{k},x_{\sigma(k)})\,\dd x_k = \tilde{K}_{n}(x_{\sigma^{-1}(k)},x_{\sigma(k)}).$$ This reconstructs $\prod_{i = 1}^{k - 1} \tilde{K}_{n}(x_i,x_{\tilde{\sigma}(i)})$, which only depends on the permutation $\tilde{\sigma} \in \mathfrak{S}_{k - 1}$ obtained from $\sigma$ by ``jumping over $k$'', i.e. $\tilde{\sigma}(i) = \sigma(i)$ if $i \neq \sigma^{-1}(k)$, and $\tilde{\sigma}(\sigma^{-1}(k)) = \sigma(k)$. There are exactly $(k - 1)$ ways to obtain a given $\tilde{\sigma}$ from some $\sigma$, since we have to choose the position of the element $\sigma(k) \in \{1,\ldots,k - 1\}$. Besides, we have ${\rm sgn}(\tilde{\sigma}) = -{\rm sgn}(\sigma)$ since the length of one cycle in $\tilde{\sigma}$ was reduced by $1$ compared to $\sigma$. All in all, these terms reconstruct:
$$
- (k - 1)\,\mathop{{\rm det}}_{1 \leq i,j \leq k - 1}\big[\tilde{K}_n(x_i,x_j)\big]\,.
$$
\end{itemize}
Summing the two entails the claim. \hfill $\Box$

\subsubsection{Spectral density}

The formula \eqref{DetP} is remarkable: we say that the eigenvalues of hermitian matrices in invariant ensembles form a \textbf{determinantal point process}. If $\tilde{K}_{n}$ was an arbitrary function of two variables, the $k \times k$ determinant of $\tilde{K}_n(x_i,x_j)$ would have no reason to be non-negative. Here, for the Christoffel-Darboux kernel, it must be non-negative by consistency.

For instance, the exact spectral density is $1/n$ times
\beq
\label{spde}\rho_{n}^{(1)}(x) = \tilde{K}_n(x,x) = \frac{p_{n}'(x)p_{n - 1}(x) - p_{n - 1}'(x)p_{n}(x)}{h_{n - 1}}\,e^{-W(x)}\,.
\eeq

\subsubsection{In the GUE}

The GUE corresponds to the weight:
\beq
\label{WG} W(x) = \frac{Nx^2}{2\sigma^2},\qquad{\rm with}\,\, N = n\,.
\eeq
We have written $N$ here instead of $n$, to stress that the size of the matrix appears in two places: first, in the orthogonality weight since $W$ depends on $N = n$, and then in the degree $n$ or $(n - 1)$ of the orthogonal polynomials we need to use in \eqref{CKD}. To avoid confusion, we may just perform all computations with $N$, and at the end set $N = n$ to retrieve the GUE normalized as in Section~\ref{GaussianS}. We will also choose $\sigma = 1$.

The orthogonal polynomials for the weight $\dd x\,e^{-x^2/2}$ on $\mathbb{R}$ are well-known, called the \textbf{Hermite polynomials} and denoted $H_n(x)$. The orthogonal polynomials for the weight $\dd x\,e^{-W(x)}$ with \eqref{WG} are just:
\beq
\label{PnPX} P_n(x) = N^{-n/2}\,H_n(N^{1/2}x)\,.
\eeq
We list basic properties of the Hermite polynomials, that can be easily derived using the orthogonality relations:
\begin{itemize}
\item[$\bullet$] $H_n$ has parity $(-1)^n$.
\item[$\bullet$] We have the formula $H_n(x) = (-1)^n\,e^{x^2/2}\,\partial_{x}^n(e^{-x^2/2})$.
\item[$\bullet$] $H_n'(x) = nH_{n - 1}(x)$.
\item[$\bullet$] We have the three-term recurrence relation $H_{n + 1}(x) = xH_{n}(x) - nH_{n - 1}(x)$.
\item[$\bullet$] The norm of $P_n$ given by \eqref{PnPX} is $h_n = \sqrt{2\pi}\,n!\,N^{-(n + 1/2)}$.
\end{itemize}

\begin{figure}[h!]
\begin{center}
\includegraphics[width=1.2\textwidth]{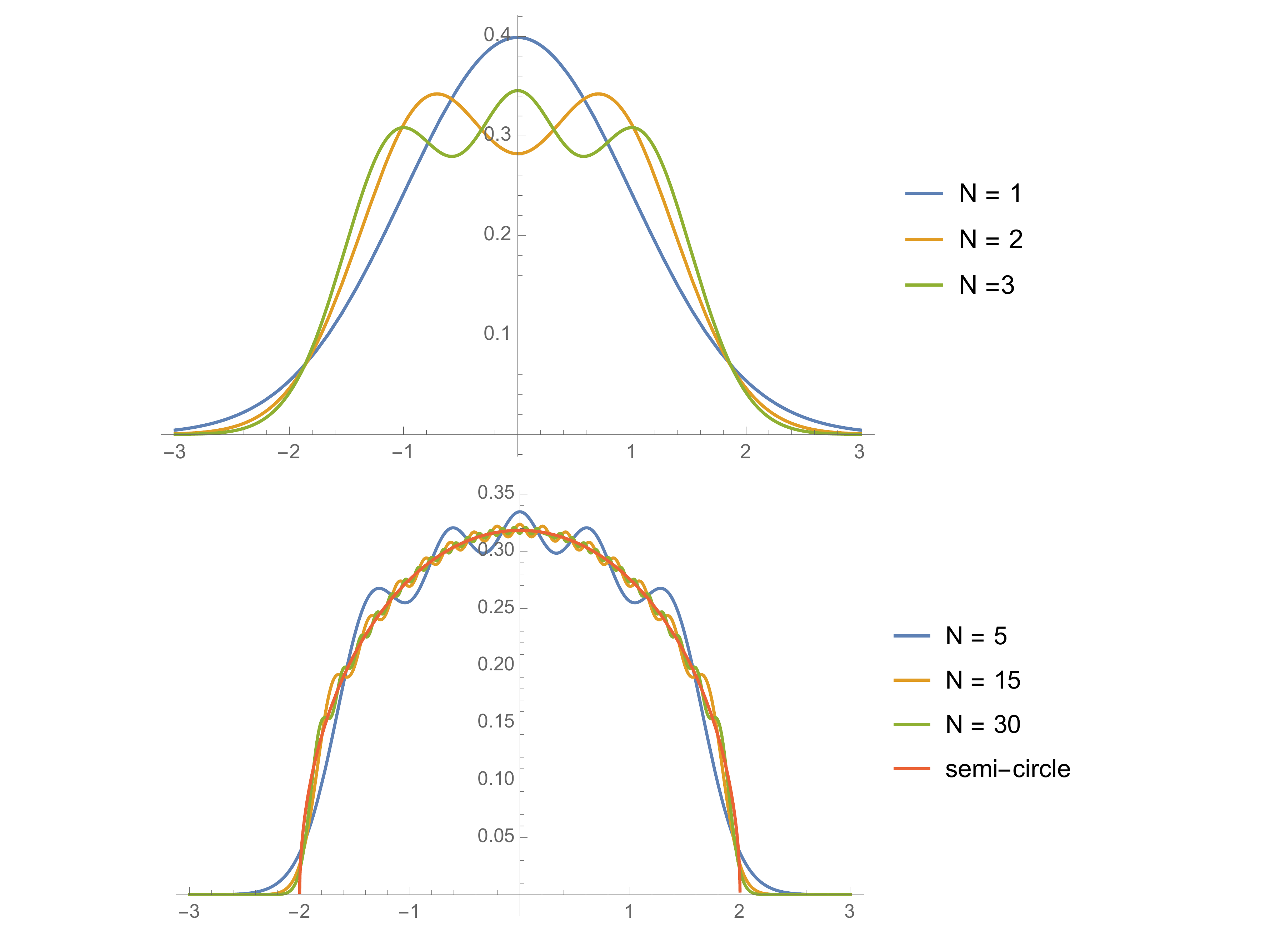}
\caption{\label{Hermite} Exact spectral density $n^{-1}\rho_{n}^{(1)}(x)$ for matrices of small size $n$ drawn from the GUE with $\sigma = 1$. For $n = 1$, this is just the Gaussian density. For $n \geq 2$ increasing, we see that it approaches the semi-circle law, with oscillations at scale $1/n$. The oscillations for $n$ finite but large can be understood as a consequence of the repulsion of eigenvalues: a region where many eigenvalues are expected prefers having less crowded neighboring regions.}
\end{center}
\end{figure}

Thus, the formula \eqref{spde} for the spectral density specializes to $1/n$ times (Figure~\ref{Hermite}):
\beq
\label{RhoGU}\rho_{n}^{(1),{\rm GUE}}(x) = \frac{\sqrt{n}}{n!\sqrt{2\pi}}\big[H_{n -1}(\sqrt{n}x)\big]^2\,e^{-nx^2/2}\,.
\eeq

\label{Exact}

\subsubsection{In the complex Wishart ensemble}

For the Wishart ensemble, one should choose an orthogonality weight on the real positive axis $\dd x\,e^{-pV(x)}$ with $V$ given by \eqref{WisV} -- and we remind that the size now is denoted $p$ instead of $n$. The corresponding orthogonal polynomials are also well-known, and called the \textbf{Laguerre polynomials}. This makes the computations in the complex Wishart ensemble rather explicit, and amenable to large $n$ asymptotics.

\subsection{Gap probabilities}

The probability that none of the eigenvalues fall into a given measurable set $A$ is also computable in terms of J\'anossy densities:
\bea
\mathbb{P}\big[{\rm no}\,\,{\rm eigenvalue}\,\,{\rm in}\,\,A\big] & = & \mathbb{E}\Big[\prod_{i = 1}^n (1 - \mathbf{1}_{A}(\lambda_i))\Big] \nonumber \\
& = & \sum_{k = 0}^{n} (-1)^k \sum_{1 \leq i_1 < \ldots < i_k \leq n} \mathbb{P}\big[\lambda_{i_1},\ldots,\lambda_{i_k} \in A\big] \nonumber \\
& = & \sum_{k = 0}^n \frac{(-1)^k}{k!} \int_{A^k} \rho_{n}^{(k)}(x_1,\ldots,x_k)\,\prod_{i = 1}^k \dd x_i\,.
\eea
From \eqref{DetP}, we find:
$$
\mathbb{P}\big[{\rm no}\,\,{\rm eigenvalues}\,\,{\rm in}\,\,A\big] = \sum_{k = 0}^{n} \frac{(-1)^k}{k!} \int_{A^k} \mathop{{\rm det}}_{1 \leq i,j \leq k}\big[\tilde{K}_{n}(x_i,x_j)\big] \prod_{i = 1}^k \dd x_i\,.
$$
Since $K_n$ is the kernel of an operator of rank $n$, the determinants of size $k > n$ vanish, and we have:
$$
\mathbb{P}\big[{\rm no}\,\,{\rm eigenvalues}\,\,{\rm in}\,\,A\big] = \sum_{k = 0}^{\infty} \frac{(-1)^k}{k!} \int_{A^k} \mathop{{\rm det}}_{1 \leq i,j \leq k}\big[\tilde{K}_{n}(x_i,x_j)\big] \prod_{i = 1}^k \dd x_i\,.
$$
We recognize the definition of the Fredholm determinant\footnote{This is a generalization of the notion of determinant to operators in infinite-dimensional spaces.} of the operator $\widehat{K}_{n}$ restricted to act on the Hilbert space $L^2(A,\dd x)$:
\beq
\label{mimum}\mathbb{P}\big[{\rm no}\,\,{\rm eigenvalues}\,\,{\rm in}\,\,A\big] = \Det\big[1 - \widehat{K}_{n}\big]_{L^2(A,\dd x)}\,.
\eeq
The Fredholm determinant $\Det[1 - \widehat{K}]$ is a continuous function of $\widehat{K}$ for the topology induced by the sup-norm for the kernel $\widehat{K}(x,y)$ of $\widehat{K}$. This means that, to study the large $n$ asymptotics of \eqref{mimum}, it is enough to study the uniform convergence of the kernel $\tilde{K}_{n}(x,y)$.

In particular, if we take $A$ to be the semi-infinite interval $(a,+\infty)$, the probability that no eigenvalue belongs to $A$ is exactly the probability that the maximum eigenvalue is smaller than $a$:
$$
\mathbb{P}[\lambda_{\max} \leq a] = \Det\big[1 - \widehat{K}_{n}\big]_{L^2\big((a,+\infty),\dd x\big)}\,.
$$

\newpage

\section{Asymptotics and universality of local regime}
\label{Asym}

We have expressed the J\'anossy densities and the gap probabilities in terms of the Christoffel-Darboux kernel:
\beq
\label{CD2}\tilde{K}_{n}(x,y) = \frac{P_{n}(x)P_{n - 1}(y) - P_{n - 1}(x)P_{n}(y)}{h_{n - 1}(x - y)}\,.
\eeq
In order to study the large $n$ limit of the eigenvalue distributions, we just need to derive the asymptotics of the orthogonal polynomials $P_{n}(x)$.

\subsection{Asymptotics of Hermite polynomials}

For Hermite polynomials, one can easily establish, from the properties previously mentioned, the integral representation:
$$
H_n(x) = {\rm i}^n\,e^{x^2/2}\int_{\mathbb{R}} \dd \zeta\,\zeta^n\,e^{-\zeta^2/2 - {\rm i}x\zeta}\,.
$$
The asymptotics of $H_n(x)$ can then be derived using the classical method of steepest descent analysis\footnote{This is a generalization in complex analysis of the Laplace method in real analysis to study the $\epsilon \rightarrow 0$ behavior integrals of the form $\int_{\mathbb{R}} e^{-f(x)/\epsilon})\dd x$.}  -- see e.g. \cite{AGZbook} for details. The result is called the Plancherel-Rotach formula -- see e.g. \cite{Szego}. Let us define:
$$
\varphi_n(x) = \frac{e^{-x^2/4}\,H_n(x)}{\sqrt{\sqrt{2\pi}\,n!}}\,.
$$

\begin{theorem}
\label{PRth}
Let $m$ be a fixed integer, and consider $n \rightarrow \infty$.
\begin{itemize}
\item[$\bullet$] \textbf{Bulk.} For fixed $x_0 \in (-2,2)$ and $X \in \mathbb{R}$, we have:
\beq
\label{uniB}\varphi_{n + m}(n^{1/2}x_0 + n^{-1/2}X) = \frac{2\cos\big[\theta_n(x_0,X,m)\big]}{n^{1/4}\sqrt{2\pi}(4 - x_0^2)^{1/4}} + O(n^{-3/4})\,,
\eeq
with:
\bea
\theta_n(x_0,X,m) & = & (n + m + 1){\rm arcsin}(x_0/2) - \frac{\pi(n + m)}{2} \nonumber \\
& & + \frac{nx_0\sqrt{4 - x_0^2}}{4} + \frac{X\sqrt{4 - x_0^2}}{2}\,. \nonumber
\eea
The result is uniform for $X$ in any compact of $\mathbb{R}$.
\item[$\bullet$] \textbf{Edge.} For fixed $X \in \mathbb{R}$, we have:
\beq
\label{uniA}\varphi_{n + m}(2n^{1/2} + n^{-1/6}X) = n^{-1/12}\,{\rm Ai}(X) + O(n^{-5/12})\,,
\eeq
where  ${\rm Ai}$ is the Airy function, i.e. the unique solution to ${\rm Ai}''(X) = X{\rm Ai}(X)$ which decays\footnote{At $X \rightarrow -\infty$, ${\rm Ai}(X)$ is unbounded and has oscillatory asymptotics.} when $X \rightarrow +\infty$ like:
$$
{\rm Ai}(X) \sim \frac{\exp\big(-\frac{2}{3}\,X^{3/2}\big)}{2\sqrt{\pi}\,X^{1/4}}\,.
$$
\eqref{uniA} is uniform for $X$ in any compact of $\mathbb{R}\cup\{+\infty\}$.
\item[$\bullet$] \textbf{Far side.} For fixed $|x_0| > 2$, $\varphi_{n + m}(n^{1/2}x_0)$ decays exponentially fast when $n \rightarrow \infty$.
\end{itemize}
\end{theorem}
The existence of the three regimes has direct qualitative consequences for the distribution of eigenvalues in the large $n$ limit. In the bulk, the Hermite polynomials have an oscillatory asymptotics: it is the region where their $n$ zeroes accumulate, and where the eigenvalue distribution will be concentrated. As expected, with the scaling \eqref{PnPX}, we look at arguments of the Hermite polynomials at the scale $\sqrt{n}$, and the bulk thus correspond to the bounded interval $x_0 \in (-2,2)$. In \eqref{uniB}, we see that non-trivial variations occur when we deviate from $x_0$ with order of magnitude $1/n$, as measured by $X$. This means that fluctuations of eigenvalues in the bulk of the GUE will occur at scale $O(1/n)$. The result in the far side indicates that it will be exponentially unlikely to find eigenvalues outside of $[-2,2]$, and confirms that the support of the spectral density should be $[-2,2]$. At the right edge $x_0 = 2$ between the far side and the bulk -- the behavior at the left edge $x_0 = -2$ is obtained by symmetry -- there is a transition, and non-trivial variations now occur when $x_0$ deviates from $2$ with order of magnitude $n^{-1/2}\cdot n^{-1/6} = n^{-2/3}$. So, the fluctuations of eigenvalues near the edge, and in particular the fluctuations of the maximum, will be of order $n^{-2/3}$, as anticipated in Section~\ref{Mxun}.

Notice that the introduction of the variable $X$ in Theorem~\ref{PRth} allows to reach the distribution of eigenvalues in regions where only finitely many eigenvalues are expected -- these are regions of size $1/n$ in the bulk, and of size $n^{-2/3}$ around the edge. I.e. it makes possible to access the local regime, while keeping only $x_0$ would provide information about the global regime only.

There is no difficulty in computing the asymptotics of the Christoffel-Darboux kernel \eqref{CD2} in the various regimes from Theorem~\ref{PRth}, although the algebra is a bit lengthy. We now summarize the results of these computations.

\subsection{Consequences in the bulk}

First, we find that the spectral density converges to the semi-circle law:
$$
\lim_{n \rightarrow \infty} n^{-1}\rho_{n}^{(1)}(x_0) = \frac{\sqrt{4 - x_0^2}}{2\pi}\,\mathbf{1}_{[-2,2]}(x_0)\,.
$$
For the local regime around a point $x_0 \in (-2,2)$ in the bulk, we find:
\beq
\label{labla} \lim_{n \rightarrow \infty} \frac{\tilde{K}_{n}\bigg(x_0 + \frac{X}{\rho_{n}^{(1)}(x_0)},x_0 + \frac{Y}{\rho_{n}^{(1)}(x_0)}\bigg)}{\rho_{n}^{(1)}(x_0)} = \frac{\sin\pi(X - Y)}{\pi(X - Y)}\,.
\eeq
This function is called the \textbf{sine kernel}, and denoted $K_{{\rm sin}}(X,Y)$. The corresponding operator is denoted $\widehat{K}_{{\rm sin}}$. In \eqref{labla}, It was natural, instead of choosing to measure $X$ in units of $1/n$, to normalize it further by the spectral density. Indeed, the average local density of eigenvalues measured in terms of $X$ is equal to $1$, and this facilitates the comparison between different models.

\begin{corollary}
For any fixed integer $k$, and fixed $x_0 \in (-2,2)$, the eigenvalue distribution is such that:
$$
\lim_{n \rightarrow \infty} \frac{\rho_{n}^{(k)}\Big[\Big(x_0 + \frac{X_i}{\rho_{n}^{(1)}(x_0)}\Big)_{i = 1}^n\Big]}{\rho_{n}^{(1)}(x_0)^k} = \mathop{{\rm det}}_{1 \leq i,j \leq k} K_{{\rm sin}}(X_i,X_j)\,.
$$
And, for any compact $A$ of $\mathbb{R}$, the gap probability behaves like:
$$
\lim_{n \rightarrow \infty} \mathbb{P}\bigg[{\rm no}\,\,{\rm eigenvalue}\,\,{\rm in}\,\,\bigg(x_0 + \frac{A}{\rho_{n}^{(1)}(x_0)}\bigg)\bigg] = \Det\big[1 - \widehat{K}_{{\rm sin}}\big]_{L^2(A,\dd x)}\,,
$$
where $a + b\cdot A$ the image of $A$ by the map $x \mapsto a + bx$.
\end{corollary}

\subsection{Consequences at the edge}

\label{Asymmax}

We find that the Christoffel-Darboux kernel at the edge behaves like:
$$
\lim_{n \rightarrow \infty} n^{-1/6}\,\tilde{K}_{n}(2 + n^{-2/3}X,2 + n^{-2/3}Y) = \frac{{\rm Ai}(X){\rm Ai}'(Y) - {\rm Ai}'(X){\rm Ai}(Y)}{X - Y}\,.
$$
This is the \textbf{Airy kernel}, denoted $K_{{\rm Ai}}(X,Y)$. The corresponding operator is denoted $\widehat{K}_{{\rm Ai}}$.

\begin{corollary}
At the right edge of the spectrum, the eigenvalue distribution is such that:
$$
\lim_{n \rightarrow \infty} n^{-k/6}\,\rho_{n}^{(k)}\big[(2 + n^{-2/3}X_i)_{i = 1}^n\big] = \mathop{{\rm det}}_{1 \leq i,j \leq k} K_{{\rm Ai}}(X_i,X_j)\,.
$$
And, for any compact $A$ of $\mathbb{R}\cup\{+\infty\}$, the gap probability behaves like:
$$
\lim_{n \rightarrow \infty} \mathbb{P}\big[{\rm no}\,\,{\rm eigenvalue}\,\,{\rm in}\,\,2 + n^{-2/3}A\big] = \Det[1 - \widehat{K}_{{\rm Ai}}]_{L^2(A,\dd x)}\,.
$$
In particular:
$$
\lim_{n \rightarrow \infty} \mathbb{P}[\lambda_{\max} \leq 2 + n^{-2/3}s\big] = \Det[1 - \widehat{K}_{{\rm Ai}}]_{L^2\big((s,+\infty),\dd x\big)}\,.
$$
is another expression -- the first historically obtained -- of the Tracy-Widom law ${\rm TW}_{2}(s)$.
\end{corollary}

\subsection{Universality}

Here is a table summarizing the  limit distributions we have encountered.
\begin{figure}[h!]
\begin{center}
\includegraphics[width=0.75\textwidth]{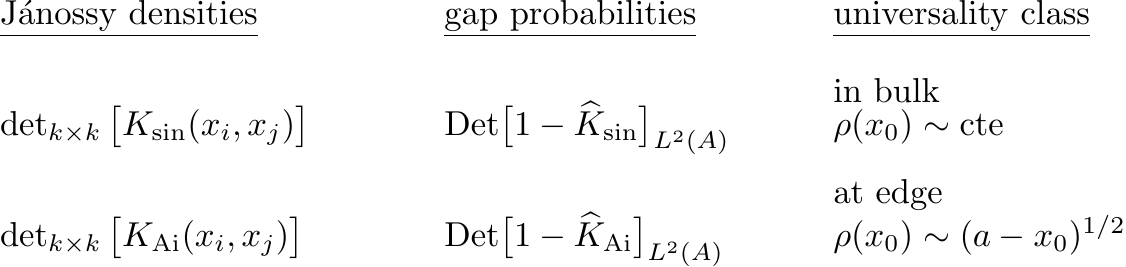}
\end{center}
\end{figure}

They are universal -- i.e. valid independently of the details of the model -- for hermitian random matrices in invariant ensembles, for complex Wishart matrices, and many other ensembles of random hermitian matrices. For symmetric matrices, there exist different universal laws -- we have seen an expression of ${\rm TW}_{1}(s)$ in \eqref{TW1}-- which are also well understood \cite{MehtaBook}. Actually, this universality goes beyond random matrices, see e.g. the review \cite{Deiftuni}. Let us illustrate it by two examples.

\subsubsection{Non-intersecting random walks}

Consider the standard brownian motion (BM) in $\mathbb{R}$, and let $\mathcal{K}_{t}(x,y)$ be the probability density that a BM starting at time $t = 0$ at position $x$, ends at time $t$ at position $y$. It is a basic result of stochastic processes that:
$$
\mathcal{K}_{t}(x,y) = (2\pi t)^{-1/2}\,\exp\Big(-\frac{(x - y)^2}{2t}\Big)\,.
$$
Since BM is a Markov process, we also have:
$$
\int_{\mathbb{R}} \mathcal{K}_{t}(x,z)\,\mathcal{K}_{t'}(z,y)\,\dd z = \mathcal{K}_{t + t'}(x,y)\,.
$$
Now, let us consider $n$ independent BMs starting from positions $x_1 < \ldots < x_n$ at time $t = 0$, which we condition not to intersect. Karlin and McGregor in 1960 \cite{KarlinMcG} have computed the probability density that they arrive at time $t$ at positions $y_1 < \ldots < y_n$:
$$
\mathcal{P}_{n}(x_1,\ldots,x_n|y_1,\ldots,y_n) = \mathop{{\rm det}}_{1 \leq i,j \leq n} \mathcal{K}_{t}(x_i,y_j)
$$

\begin{figure}[h!]
\begin{center}
\includegraphics[width=0.7\textwidth]{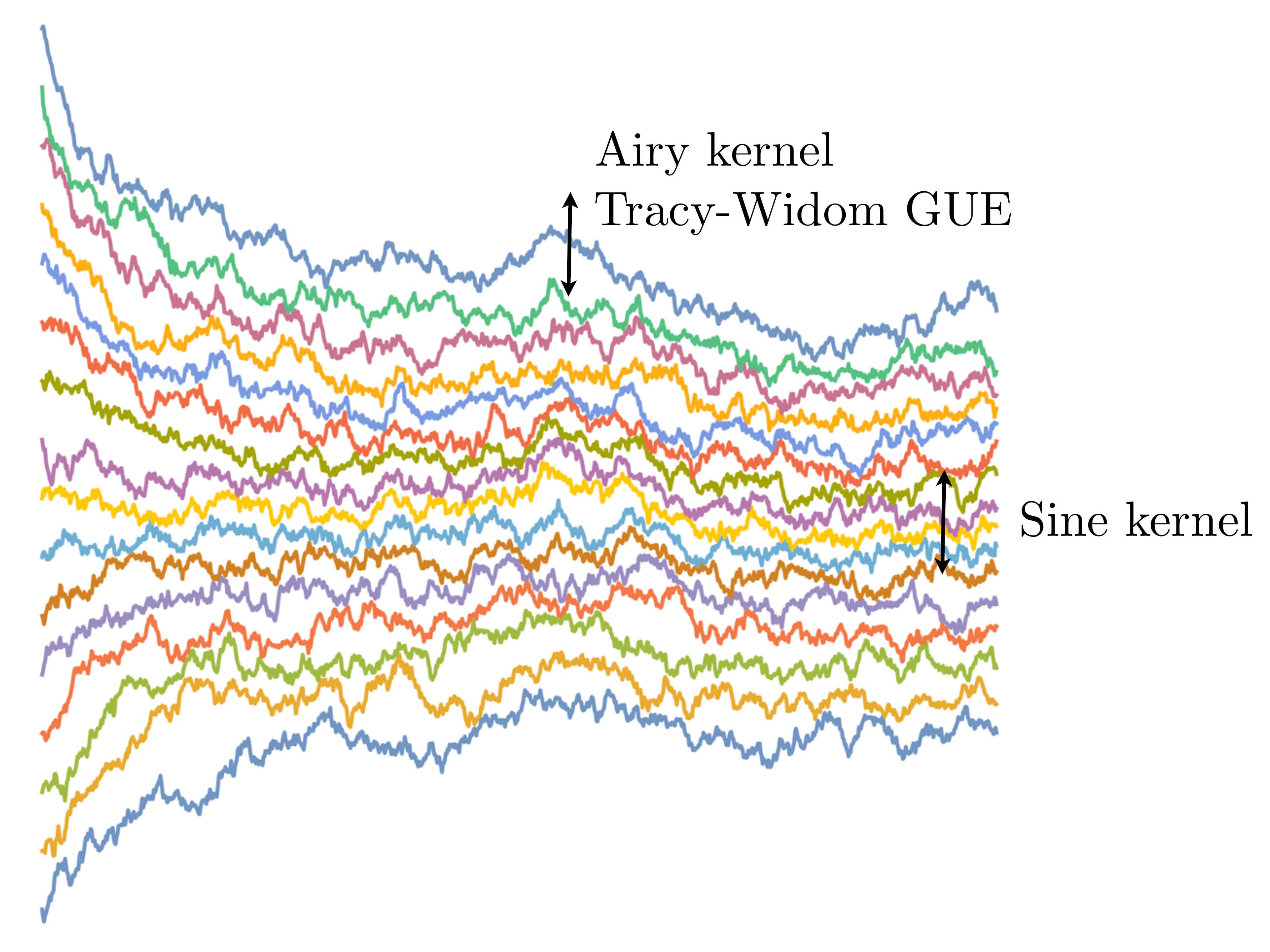}
\caption{\label{Nonint} Simulation (courtesy of P.~Ferrari) $n$ independent random walks in dimension one, conditioned not to intersect. In the large $n$ limit, after proper rescaling, the fluctuations of the height of the top path follows the Tracy-Widom GUE law, and the joint distribution of a finite number of paths starting from the top path is given by the determinantal process with kernel $K_{{\rm Ai}}$. For a path in the bulk, the fluctuations of the height of a finite number of consecutive paths are given by the  determinantal process with kernel $K_{{\rm sin}}$.}
\end{center}
\end{figure}

This is the starting point of a series of results, showing that in various situations, the non-intersecting random walkers -- sometimes called vicious because they do not want to cross -- behave when $n \rightarrow \infty$ like eigenvalues of large random matrices (Figure~\ref{Nonint}). For instance, the fluctuations of the position of the rightmost walker generically occur at scale $n^{-2/3}$ around their mean, and converge in law towards the Tracy-Widom GUE law. Similarly, if one zooms amidst the walkers in a region where we expect to see only finitely many of them, the distribution of the positions of $k$ of them is given by the $k \times k$ determinant built from the sine kernel. More details can be found in \cite{Ferrarilec}.

\subsubsection{Growth models}

The sine kernel or the Airy kernel distributions also appear in problems of growing interfaces. There exist several mathematical models where this has been established -- see the review \cite{FerrariSpohn}. But I also want to point out, with an example, that these distributions can be seen in (even non-mathematical) nature.

\begin{figure}[h!]
\begin{center}
\includegraphics[width=\textwidth]{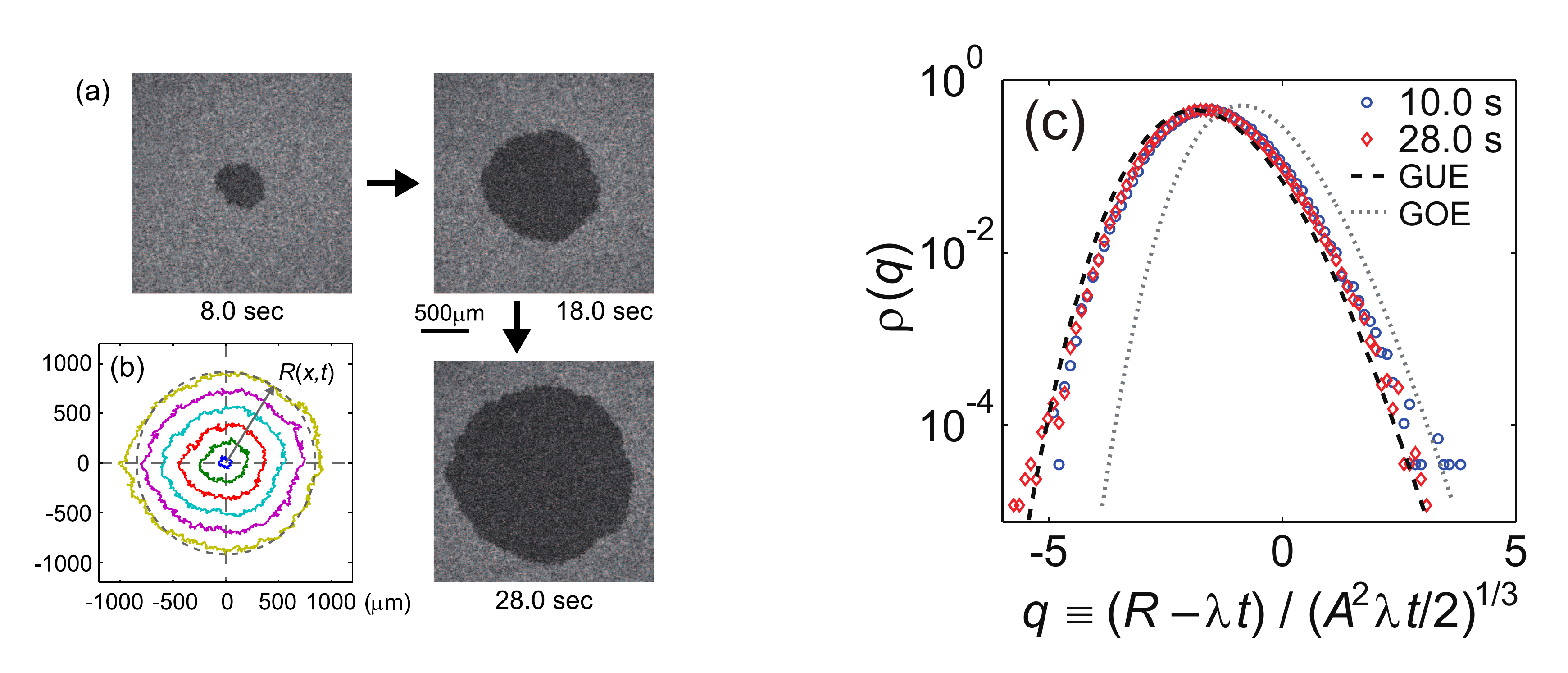} 
\caption{\label{TSano} Comparison between fluctuations of the radius of a growing interface in nematic liquid crystals and Tracy-Widom laws. Reprinted with permission from \emph{Universal fluctuations of growing interfaces: evidence in turbulent liquid crystals}, K.~Takeuchi and M.~Sano, Phys. Rev. Lett. \textbf{104} 230601 (2010) \copyright\,APS.}
\end{center}
\end{figure}

The physicists Takeuchi and Sano (2010) observed experimentally the Tracy-Widom law in nematic liquid crystals. ``Nematic'' means that the material is made of long molecules whose orientation has long-range correlations, while liquid means that the molecules in the neighborhood of a given one are always changing, i.e. the correlation of positions have short range. In nematic materials, a ``topological defect'' is a configuration of orientations that winds around a point. In two dimensions, it occurs for instance when the local orientation rotates like the tangent vector when following a circle throughout the material\footnote{In three dimensions, the Hopf fibration $\phi\,:\,\mathbb{S}_{3} \rightarrow \mathbb{S}_2$ is a configuration of orientations realizing a topological defect.}. The material studied by Takeuchi and Sano admits two phases: the phase appearing here in gray (resp. black) has a low (resp. high) density of topological defects. If one applies a voltage to the grey phase, one encourages the formation of defects. Once this happens -- here at the center of the picture at time $t = 0$ -- the black phase takes over the grey phase from this primary cluster of defects. One observes that the interface grows approximately linearly with time $t$. However, the turbulence driving the system causes some fluctuations from samples to samples. The distribution of these fluctuations of radius around the linear drift matches with the Tracy-Widom GUE law, and the quality of the fit improves with time increasing (Figure~\ref{TSano}). The symmetry class in this case is conditioned by the geometry: a spherical geometry leads to GUE, while a flat interface between two phases would lead to GOE. This result is confirmed in a mathematical model for the interface growth analyzed at $t \rightarrow +\infty$ by Sasamoto and Spohn around the same time \cite{SpohnKPZ}.

\subsubsection{Last remarks}

In the last twenty years, tremendous progress has been made to prove universality in random matrices, with weak assumptions, relying on various approaches. Without exhaustivity, we can cite:
\begin{itemize}
\item[$\bullet$] the fact that some models are exactly solvable (like the invariant ensembles of symmetric or hermitian random matrices) and Riemann-Hilbert steepest descent analysis. This is very useful, but maybe not very satisfactory from the probabilistic point of view, since the method hinges from the beginning on ``algebraic miracles'', which are not anymore available if the models are slightly perturbed.
\item[$\bullet$] transport of measures (Shcherbina ; Figalli, Guionnet and Bekerman), which has succeeded in proving some universality for all $\beta$-ensembles.
\item[$\bullet$] relaxation methods (Bourgade, Erd\"os, H.-T. Yau, etc.) which are purely based on probability, stochastic processes and analysis, and brought many results for invariant ensembles, matrices with independent entries, etc.
\item[$\bullet$] combinatorial methods (Wigner ; Soshnikov ; Tao and Vu, etc.) which are particularly useful for matrices with independent entries, etc.
\end{itemize}
One current trend is now to apply the insight gained from the study of random matrices, to more difficult problems like random band matrices, random Schr\"odinger operators, adjacency matrices of random graphs, etc. This is motivated by the desire to understand the properties of localization/delocalization of the eigenvectors -- that determine isolating/conducting properties of materials modelized in this way.

\newpage

\section{Questions of participants}

\noindent $\bullet$ \textbf{Ninjbat Uuganbaatar: Can one apply PCA techniques to analyze voting?}

In general, the number of options for which one can vote is very small, so I do not see how PCA can be used to analyze voting. However, it could be a tool to check the representativity of the political offer in a given society. For instance, one could ask $n$ individuals to answer a poll consisting of $p$ questions about their political preferences. As example of questions: how much should income be taxed? at which age should people retire? should the state subsidize health coverage? \ldots{} The opinion pollster would have to choose a way to get answers which are numbers, for instance binary questions -- somewhat like in population genetics about presence or absence of an allele -- given $0$ or $1$ as entries, or questions that one can answer by an intensity from $0$ (not at all) to $10$ (absolutely). Then, one can build a $n \times p$ matrix $X$ collecting the answers, and the empirical covariance matrix $M = p^{-1}\,XX^{T}$. By PCA analysis, one can then hope to determine how many relevant groups can be formed, that have similar political ideas -- as probed by the questions asked. One could then compare with the number of political parties, as well as their programme, to see if the population is well-represented at the level of ideas, and if their strength compares well with the magnitude of the eigenvalues found in PCA. I do not know if such a project has been already conducted. Clearly, an important work of calibration  is needed -- e.g. checking if the outcome of PCA is similar when one asks yes/no questions, or intensity questions, etc. -- to ensure the results are reliable.

\vspace{0.3cm}

\noindent $\bullet$ \textbf{Remco van der Hofstad: How can one identify quantitatively in PCA what comes from true information and what comes from noise?}

\vspace{0.2cm}

For market prices, we have seen in the examples of Section~\ref{MarS} that the overlap between the $j$-th eigenvector -- sorted by decreasing order for the corresponding eigenvalues -- of empirical correlation matrices in two distinct periods does not exceed what one expects from the overlap of two independent random vectors for some $j \geq j_0$. And this threshold also corresponded well with the position of the noise band -- i.e. the distribution of eigenvalues $\lambda_j$ with $j \geq j_0$ was fit with the Mar\v{c}enko-Pastur law.

A more general method is to fix a confidence threshold, and then make a statistical test for $\lambda_i$ using the Tracy-Widom law, for $i = 1,2,3,\ldots$ until one cannot reject anymore the null hypothesis (which enjoys Tracy-Widom distribution). More precisely, if the test is passed for $\lambda_i$, one restricts the matrix to the orthogonal of the eigenspace of $\lambda_{1},\ldots,\lambda_i$ before continuing the analysis. And there exists estimates of the rate of convergence to the Tracy-Widom law in null Wishart matrices (see e.g. \cite{JohnTW}) when $n,p$ is large but not infinite, which can be used for statistical tests. To cope with finite size effects, one can also use large deviation functions -- see the question below -- but one should keep in mind that their details are much less robust (if one changes the model) than the Tracy-Widom distribution.

\newpage
\vspace{0.3cm}

\noindent $\bullet$ \textbf{Kanstantsin Matetski: What can be said about the large deviations of the maximum eigenvalue?}

\vspace{0.2cm}

Although I did not present them for lack of space in the lectures, there exist techniques, based on potential theory and large deviation theory, to compute the asymptotic behavior of the partition function in invariant ensembles. In particular, if one assume that the support of the large $n$ spectral density is a single segment (as for GUE and Wishart) + some other technical assumptions on $V$, one can show that the partition function:
$$
Z_{n,\beta}(A) = \int_{A^n} \prod_{1 \leq i < j \leq n} |\lambda_j - \lambda_i|^{\beta}\,\prod_{i = 1}^n \exp\Big(-\frac{n\beta}{2}\,V(\lambda_i)\Big)\dd\lambda_i
$$
has an asymptotic expansion of the form:
\beq
\label{asu}\ln Z_{n,\beta}(A) = n^2 F_0 + (\beta/2)n \ln n + n(\beta/2 - 1) F_1 + \frac{3 + 2/\beta + \beta/2}{12}\ln n + F_2 + o(1)
\eeq
when $n \rightarrow \infty$, and the coefficients $F_j$ can be computed fairly explicitly, depending on $V$ and $A$. The $o(1)$ actually consists of a full asymptotic expansion in powers of $1/n$, and its coefficients can also be computed recursively.

These results give access to the large deviations for the maximum eigenvalue, since:
$$
\mathbb{P}[\lambda_{\max} \leq a] = \frac{Z_{n,\beta}(a,+\infty)}{Z_{n,\beta}(\mathbb{R})}\,.
$$
For instance, when $a$ is independent of $n$ and strictly smaller than $a_* = \lim_{n \rightarrow \infty} \mathbb{E}[\lambda_{\max}]$, the assumptions leading to \eqref{asu} are satisfied and we can prove rigorously an asymptotic expansion of the form:
\beq
\label{LL}\mathbb{P}[\lambda_{\max} \leq a] = n^{c}\,\exp\Big[-n^2G_0(a) - n(\beta/2 - 1)G_1(a) - \sum_{k \geq 0}^{K} n^{-k}\,G_{k + 2}(a) + o(n^{-K})\Big]\,.
\eeq
For $a < a_*$, this probability is super-exponentially small because one has to push all the $n$ eigenvalues to the left of $a_*$ to achieve the event $\lambda_{\max} \leq a < a_*$. The leading term $G_0(a)$ is called the large deviation function, and has some relevance in statistical applications, because one has to face the finite size of data.

How does that connect to the Tracy-Widom law? If one naively inserts $a = a_* - sn^{-2/3}$ in the right-hand side of \eqref{LL}, we can show that each term $n^{-k}\,G_{k + 2}(a)$ tends to a constant $\tilde{G}_{k + 2}\,s^{-3k/2}$, which is of order $1$. This is not surprising because in this regime the probability \eqref{LL} should vary between $0$ and $1$. As a matter of fact, putting $a = a_* - sn^{-2/3}$ goes out of the range in which \eqref{LL} was established. But, if one is ready to believe that the crossover from ``large deviations'' to ``not so large deviations'' is smooth -- an exchange of limits that has not been justified as of writing -- then we interpret the naive right-hand side where one first inserts $a = a^* - sn^{-2/3}$ as the all-order asymptotic expansion when $s \rightarrow +\infty$ of ${\rm TW}_{\beta}(-s)$. This leads to predictions, for any value of $\beta > 0$, for the left tail of Tracy-Widom $\beta$ laws. They agree with all rigorous results known for $\beta = 1,2$, and with the leading order rigorously known for arbitrary $\beta$. In particular, we have a prediction for the constant term of the asymptotic expansion, which is always tricky to get. A similar story can be devised for the right tail.

The large deviation function $G_0(a)$ at the left tail was first computed by Dean and Majumdar in \cite{DeanMaj} -- although this is a physics paper, the equation they solve to get $G_0(a)$ can be rigorously established using potential theory without any difficulty, hence making a complete proof. We discussed the generalization to all-order finite size corrections in \cite{BEMN} for the left tail, and \cite{BNbeta} for the right tail. The computations in these two papers are done for the Gaussian ensembles, but there would be no difficulty in conducting them for other $V$, e.g. for the Wishart ensembles. These two papers take as starting point the asymptotic expansion of the form \eqref{LL} ;  these expansions have been established rigorously in \cite{BG11}.

\vspace{0.3cm}

\noindent $\bullet$ \textbf{Ninjbat Uuganbaatar: Is there a combinatorial interpretation to the formulas we have seen for the distribution of random matrices?}

Let us start with a matrix $M_n$ in the Gaussian ensembles, for $\sigma = 1$. The moments of the semi-circle law can be directly computed by expanding its Stieltjes transform \eqref{WSClaw} at $z \rightarrow \infty$:
$$
\lim_{n \rightarrow \infty} n^{-1}\,\mathrm{Tr}\,M_n^{2k} = \frac{2k!}{k!(k + 1)!} = {\rm Cat}(k)\,.
$$
This is the Catalan number, computing the number of ways to connect pairs of edges in a $2k$-gon, without crossing. More generally, Harer and Zagier in 1986 \cite{HarerZagier} showed the expansion:
$$
\mathrm{Tr}\,M_n^{2k} = \sum_{g \geq 0} n^{1 - 2g}\,\mathcal{N}_{n}(g)
$$
where $\mathcal{N}_{n}(g)$ is the number of ways of identifying by pairs the edges of $2k$-gon, in such a way that the resulting surface has genus $g$. They gave several formulas to compute these numbers -- from \eqref{RhoGU}, we know that they can be expressed in terms of Hermite polynomials. Harer and Zagier used this to compute the Euler characteristics of the moduli space of Riemann surfaces of genus $g$ ; this is one of the many and fruitful point of contacts between random matrices and algebraic geometry.

Actually, the combinatorial interpretation of the moments of the GUE was already known to physicists, in the more general context of invariant ensembles of hermitian matrices. Br\'ezin, Itzykson, Parisi and Zuber showed in 1979 \cite{BIPZ} that the partition function ``decomposes'' as:
$$
Z_n = n^{n + 5/12}\,\exp\Big(\sum_{g \geq 0} n^{2 - 2g}\,\mathcal{F}_{g}\Big)\,,
$$
and $\mathcal{F}_{g}$ enumerates discretized surfaces of genus $g$. For instance, if one takes $V(x) = x^2/2 - t x^3/3$, $\mathcal{F}_{g}$ is the number of triangulations of a genus $g$ surface, counted with a weight $t^{T}$ if it is made exactly of $T$ triangles. Although it seems he partition function does not make sense as a convergent integral since $V(x) \rightarrow -\infty$ when $x \rightarrow {\rm sgn}(t)\infty$, it can be defined rigorously as a formal series in the parameter $t$ -- and this is why I said ``decompose'' with quotes. Likewise the expectation values:
$$
\mathbb{E}\big[{\rm Tr}\,M_n^{\ell_1} \cdots \mathrm{Tr}\,M_n^{\ell_k}\big]
$$
are related to the enumeration of discretized surfaces with $k$ boundaries of respective perimeters $\ell_1,\ldots,\ell_k$ counted with a weight $n^{\chi}$ where $\chi$ is the Euler characteristics. The coupling of the matrix size with the Euler characteristics is a phenomenon that was first observed in gauge theories by the theoretical physicist t'Hooft in 1974 \cite{tHooft}. More on the relations between random matrices, enumeration of discretized surfaces and algebraic geometry, can be found in the book \cite{Ebook}.

\newpage

\providecommand{\bysame}{\leavevmode\hbox to3em{\hrulefill}\thinspace}
\providecommand{\MR}{\relax\ifhmode\unskip\space\fi MR }
\providecommand{\MRhref}[2]{%
  \href{http://www.ams.org/mathscinet-getitem?mr=#1}{#2}
}
\providecommand{\href}[2]{#2}


\begin{thebibliography}{10}

\bibitem{AGZbook}
G.W. Anderson, A.~Guionnet, and O.~Zeitouni, \emph{An introduction to random
  matrices}, Cambridge University Press, 2010.

\bibitem{BBPphaset}
J.~Baik, G.~Ben Arous, and S.~P{\'e}ch{\'e}, \emph{Phase transition of the
  largest eigenvalue for nonnull complex sample covariance matrices}, Ann.
  Probab. \textbf{33} (2005), 1643--1697, math.PR/0403022.

\bibitem{VirBloe}
A.~Bloemendal and B.~Vir{\'{a}}g, \emph{Limits of spiked random matrices {I}},
  Probab. Th. Rel. Fields \textbf{156} (2013), no.~3-4, 795--825,
  math.PR/1011.1877.

\bibitem{BEMN}
G.~Borot, B.~Eynard, S.N. Majumdar, and C.~Nadal, \emph{Large deviations of the
  maximal eigenvalue of random matrices}, J. Stat. Mech. (2011), no.~P11024,
  math-ph/1009.1945.

\bibitem{BG11}
G.~Borot and A.~Guionnet, \emph{Asymptotic expansion of $\beta$ matrix models
  in the one-cut regime}, Commun. Math. Phys \textbf{317} (2013), no.~2,
  447--483, math.PR/1107.1167.

\bibitem{BNbeta}
G.~Borot and C.~Nadal, \emph{Right tail expansion of {T}racy-{W}idom beta
  laws}, RMTA \textbf{1} (2012), no.~03, math-ph/1111.2761.

\bibitem{BIPZ}
\'{E}. Br{\'{e}}zin, C.~Itzykson, G.~Parisi, and J.-B. Zuber, \emph{Planar
  diagrams}, Commun. Math. Phys. \textbf{59} (1978), 35--51.

\bibitem{DeanMaj}
D.S. Dean and S.N. Majumdar, \emph{Large deviations of extreme eigenvalues of
  random matrices}, Phys. Rev. Lett. \textbf{97} (2006), 160--201,
  cond-mat/0609651.

\bibitem{Defcours}
P.~Deift, \emph{Orthogonal polynomials and random matrices : a
  {R}iemann-{H}ilbert approach}, AMS, New York, 1998, Courant Institute of
  Mathematical Sciences.

\bibitem{Deiftuni}
\bysame, \emph{Universality for mathematical and physical systems}, Proceeding
  of the ICM, Madrid 2006, Spain (2007), 125--152, math.ph/0603038.

\bibitem{Ebook}
B.~Eynard, \emph{Counting surfaces}, Progress in Mathematics, Birkh\"auser,
  2016, available at http://eynard.bertrand.voila.net/TOCbook.htm.

\bibitem{Ferrarilec}
P.L. Ferrari, \emph{Why random matrices share universal processes with
  interacting particle systems ?},  (2013), ICTP Lecture notes,
  math-ph/1312.1126.

\bibitem{FerrariSpohn}
P.L. Ferrari and H.~Spohn, \emph{Random growth models},  (2010),
  math.PR/1003.0881.

\bibitem{WFisher}
R.A. Fisher, \emph{The sampling distribution of some statistics obtained from
  non-linear equations}, Ann. Eugenics \textbf{9} (1939), 238--249.

\bibitem{ForresterS}
P.J. Forrester, \emph{The spectrum edge of random matrix ensembles}, Nucl.
  Phys. B (1993), 709--728.

\bibitem{Geman}
S.~Geman, \emph{A limit theorem for the norm of random matrices}, Ann. Probab.
  \textbf{8} (1980), no.~2, 252--261.

\bibitem{WGirschik}
M.A. Girshick, \emph{On the sampling theory of roots of determinantal
  equations}, Ann. Math. Stat. \textbf{10} (1939), 203--204.

\bibitem{HarerZagier}
J.~Harer and D.~Zagier, \emph{The {E}uler characteristics of the moduli space
  of curves}, Invent. Math. \textbf{85} (1986), 457--485.

\bibitem{HMcL}
S.P. Hastings and J.B. Mc{L}eod, \emph{A boundary value problem associated with
  the second {P}ainlev{\'{e}} transcendent and the {K}orteweg-de {V}ries
  equation}, Archive for {R}ational {M}echanics and {A}nalysis \textbf{73}
  (1980), no.~1, 31--51.

\bibitem{Hotelling}
H.~Hotelling, \emph{Analysis of a complex of statistical variables into its
  principal components}, Journal of Educational Psychology (1931), 417--441.

\bibitem{WHsu}
P.L. Hsu, \emph{On the distribution of roots of certain determinantal
  equations}, Ann. Eugenics \textbf{9} (1939), 250--258.

\bibitem{JohnTW}
I.M. Johnstone, \emph{On the distribution of the largest eigenvalue in
  principal components analysis}, Ann. Stat. \textbf{29} (2001), no.~2,
  295--327.

\bibitem{KarlinMcG}
S.~Karlin and J.~McGregor, \emph{Coincidence probabilities}, Pacific J. Math.
  \textbf{9} (1959), no.~4, 1141--1164.

\bibitem{Markowitz}
H.~Markowitz, \emph{Portfolio selection}, J. Finance \textbf{7} (1952), no.~1,
  77--91.

\bibitem{MPbib}
V.A. Mar\v{c}enko and L.A. Pastur, \emph{Distribution of eigenvalues for some
  sets of random matrices}, Mat. Sb. \textbf{72} (1967), no.~4, 507--536.

\bibitem{MehtaBook}
M.L. Mehta, \emph{Random matrices}, 3$^{rd}$ ed., Pure and {A}pplied
  {M}athematics, vol. 142, Elsevier/Academic, Amsterdam, 2004.

\bibitem{Pearson}
K.~Pearson, \emph{On lines and planes of closest fit to systems of points in
  space}, Philosophical Magazine \textbf{2} (1901), 559--572.

\bibitem{WRoy}
S.N. Roy, \emph{{$p$}-statistics or some generalizations in the analysis of
  variance appropriate to multivariate problems}, Sankhya \textbf{4} (1939),
  381--396.

\bibitem{SpohnKPZ}
H.~Spohn and T.~Sasamoto, \emph{The one-dimensional {KPZ} equation: an exact
  solution and its universality}, Phys. Rev. Lett. \textbf{104} (2010),
  cond-mat.stat-mech/1009.1883.

\bibitem{Szego}
G.~Szeg\"o, \emph{Orthogonal polynomials}, Amer. Math. Soc., 1939, reprinted
  with corrections (2003).

\bibitem{Taobook}
T.~Tao, \emph{Topics in random matrix theory}, Graduate Studies in Mathematics,
  vol. 132, AMS, 2012.

\bibitem{tHooft}
G.~t'Hooft, \emph{A planar diagram theory for strong interactions}, Nucl. Phys.
  B \textbf{72} (1974), 461--473.

\bibitem{TW92}
C.~Tracy and H.~Widom, \emph{Level spacing distributions and the {A}iry
  kernel}, Commun. Math. Phys. \textbf{159} (1994), 151--174, hep-th/9211141.

\bibitem{TW95}
\bysame, \emph{On orthogonal and symplectic matrix ensembles}, Commun. Math.
  Phys. \textbf{177} (1996), 727--754, solv-int/9509007.

\bibitem{Voiculescu}
D.V. Voiculescu, \emph{Limit laws for random matrices and free products},
  Invent. Math. \textbf{104} (1991), 201--220.

\bibitem{Wignero}
E.P. Wigner, \emph{Characteristic vectors of bordered matrices with infinite
  dimensions}, Ann. Math. \textbf{62} (1955), no.~3, 548--564.

\bibitem{Wishart}
J.~Wishart, \emph{The generalised product moment distribution in samples from a
  normal multivariate population}, Biometrika \textbf{20A} (1928), no.~1/2,
  32--52.

\end{thebibliography}
\end{document}